\documentclass[a4paper,12pt] {article}
\usepackage{helvet}
\usepackage{bbold}
\usepackage[usenames]{color}
\usepackage{latexsym,graphics}
\usepackage[hang]{subfigure}
\usepackage{graphics}
\usepackage[dvips]{graphicx}
\usepackage{epsfig}
\usepackage{pst-grad} % For gradients
\usepackage{pst-plot} % For axes
\usepackage{pstricks}
\usepackage{lmodern}
\usepackage{xcolor,rotating,epic,eepic}
\usepackage{subfigure}
\usepackage{lscape}
\usepackage{url}
\usepackage{array}
\usepackage{multirow}
\usepackage{mathrsfs}
\usepackage{amsmath}
\usepackage{float}
\usepackage{placeins}
\usepackage{amsfonts}
\usepackage[utf8]{inputenc}
\usepackage[T1]{fontenc}
\usepackage[french, english]{babel}
\usepackage{amssymb}
\usepackage{epstopdf}

\DeclareGraphicsExtensions{.jpg, .pdf, .png, .fig, .eps}

\begin{document}

\begin{center}
\Large{\textbf{The Harmonic Balance Method for Bifurcation Analysis of Large-Scale Nonlinear Mechanical Systems}}\vspace{1cm}
\end{center}

\begin{center}
{T. Detroux, L. Renson, L. Masset, G. Kerschen\\\vspace{0.8cm}

\small Space Structures and Systems Laboratory\\
Department of Aerospace and Mechanical Engineering\\
University of Li\`{e}ge, Li\`{e}ge, Belgium \\
E-mail: tdetroux, l.renson, luc.masset, g.kerschen@ulg.ac.be\\\vspace{0.5cm} \vspace{1cm}

\rule{0.85\linewidth}{.3pt}
\vspace{-0.5cm}
\begin{abstract}

The harmonic balance (HB) method is widely used in the literature for analyzing the periodic solutions of nonlinear mechanical systems. The objective of this paper is to exploit the method for bifurcation analysis, i.e., for the detection and tracking of bifurcations of nonlinear systems. To this end, an algorithm that combines the computation of the Floquet exponents with bordering techniques is developed. A new procedure for the tracking of Neimark-Sacker bifurcations that exploits the properties of eigenvalue derivatives is also proposed. The HB method is demonstrated using numerical experiments of a spacecraft structure that possesses a nonlinear vibration isolation device.

\vspace{1cm}

\end{abstract}
\vspace{-0.5cm}
\rule{0.85\linewidth}{.3pt}

\vspace{0.5cm} Corresponding author: \\ Thibaut Detroux\\
Space Structures and Systems Laboratory\\
Department of Aerospace and Mechanical Engineering\\
University of Li\`{e}ge
\\ 1 Chemin des Chevreuils (B52/3), B-4000 Li\`{e}ge, Belgium. \\
Email: tdetroux@ulg.ac.be
\vspace{2cm}\\}
\end{center}

\normalsize

\newpage
%%%%%%%%%%%%%%%%%%%%%%%%%%%%%%%
% ---------------- INTRODUCTION
%%%%%%%%%%%%%%%%%%%%%%%%%%%%%%%
\section{Introduction}\label{intro}

Nonlinear systems are known to exhibit rich and complex dynamical behaviors, which linear systems cannot. These behaviors include, for instance, modal interactions, detached resonance curves, quasiperiodic oscillations, bifurcations and chaos. Even though periodic solutions represent only a subset of the dynamical attractors of nonlinear systems, their computation and interpretation usually provide great insight into the system's dynamics. Different algorithms and numerical methods for the computation of periodic solutions can be found in the literature. Most of them build on a continuation procedure \cite{padmanabhan1995}, as a means of studying the evolution of the periodic solutions with respect to the frequency of the harmonic forcing or a design parameter.

Time-domain methods, which deal with the resolution of a boundary value problem (BVP), have proven  effective for low-dimensional problems. When applied to larger systems, their computational burden can become substantial. For example, the shooting technique requires numerous time integrations that can slow down the algorithm. Efforts have been undertaken to make shooting less computationally intensive by using parallelization \cite{stoykov} and sensitivity analysis \cite{peeters2009}. Methods based on orthogonal collocation are utilized in several software for bifurcation detection and tracking, e.g., \textsc{auto} \cite{doedel1997}, \textsc{colsys} \cite{ascher1979}, \textsc{content} \cite{kuznetsov1995}, \textsc{matcont} \cite{dhooge2003}, and, more recently, \textsc{coco} \cite{dankowicz2011}. In spite of its high accuracy and ability to address stiff problems, orthogonal collocation is rarely employed for large systems, which can be explained by the considerable memory space required for the discretization of the BVP.

Among all methods in the frequency domain, harmonic balance (HB) is certainly the most widely used method. It is also known as the Fourier--Galerkin method, since it consists in the application of the Galerkin method with Fourier basis and test functions. The periodic signals are approximated with their Fourier coefficients, which become the new unknowns of the problem. The term `harmonic balance' was first introduced by Krylov and Bogoliubov \cite{krylov1943} who performed linearization of nonlinear dynamical equations with single-harmonic approximations. In the 1960s a demonstration of the convergence of the method for Fourier approximations truncated to several harmonics was offered by Urabe \cite{urabe1965}. The main advantage of HB is when low orders of approximation are sufficient to obtain an accurate solution, which usually holds for smooth nonlinearities. In this case, the method involves algebraic equations with less unknowns than for orthogonal collocation. The reader can refer to \cite{karkar2014} for a comparison between HB and orthogonal collocation applied to smooth and nonsmooth nonlinearities.

Several improvements and adaptations to the HB method were brought in the literature during the last couple of decades. The incremental HB applies the incremental procedure to the equations of motion before the harmonic approximation and balance (see \cite{lau1981,cheung1982} or, more recently, \cite{sze2005}). Taking advantage of the fast Fourier transform, Cameron \textit{et al.} proposed the alternating frequency-time domain (AFT) method that evaluates the nonlinear terms of the equations in the time domain where their analytical expression is known \cite{cameron1989}. Since then, many studies utilized the AFT method  \cite{cardona1998,sinou2007,jaumouille2010}. Similar developments led to the hybrid frequency-time domain method, with applications to systems with dry friction \cite{guillen1999,poudou2003}. HB was also coupled with continuation schemes, e.g., with arc-length continuation \cite{vongroll2001} or the so-called asymptotic numerical method \cite{arquier2007,cochelin2009} in the \textsc{manlab} package. In \cite{grolet2012}, a new adaptive HB method was proposed for which the number of harmonics for each degree of freedom (DOF) is automatically selected during the continuation process.

HB has enjoyed numerous applications in the literature. In electrical engineering, Kundert \textit{et al.} reported its superiority over time-domain techniques for the simulation of nonlinear circuits \cite{kundert1986}. Using a well-known variant of the method based on single-harmonic approximations, the so-called describing function (DF) method, Genesio \textit{et al.} \cite{genesio1992} provided analytical expressions for regions of chaotic behavior of Lur'e systems. Other studies were carried out, e.g., for piezoelectric inertial generators \cite{stanton2012} or DC-DC converters \cite{fang2012}.

Following a pioneer work on wing-control surface flutter \cite{shen1959}, the HB method was successfully applied to aeroelastic systems, e.g., to airfoils with freeplay \cite{liu2005} and cubic stiffness \cite{lee2005}. A comparative review of HB applied to limit cycle oscillations can be found in \cite{dimitriadis2008}. More recently, fluid dynamics problems where unsteady flows are periodic in time were also tackled. In \cite{hall2002}, Hall \textit{et al.} formulated the HB method for the Navier-Stokes equations, and applications can be found for flows in multi-stage turbomachinery \cite{gopinath2007,clark2015} and helicopter blades \cite{ekici2008}.

In mechanical engineering, Cardona \textit{et al.} developed a multi-harmonic resolution scheme for vibration analysis with nonsmooth nonlinear functions \cite{cardona1994}. The HB method was then applied to realistic examples including bladed disks \cite{petrov2003}, bolted joints \cite{jaumouille2010,suss2015}, rotor/stator contacts \cite{vongroll2001}, vibro-impact systems \cite{woo2000,peter2014}, geometrically nonlinear beams \cite{lewandowski1997} and plates \cite{ribeiro1999}. The dynamics of complete vehicles \cite{barillon2013} and full-scale aircraft \cite{senechal2014} were also examined with HB, which demonstrates its effectiveness when applied to reduced finite element models. Comparisons between experiments and the results of HB simulations gave further evidence of the accuracy of the method \cite{claeys2014,claeys2014b}.

Beyond applications to classical vibration problems, HB was exploited for extending linear modal analysis to nonlinear systems \cite{kerschen2009}. Specifically, thanks to its natural filtering property, the computation of nonlinear normal modes (NNMs) and their interactions could be substantially improved \cite{krack2013b,detroux2014}. As a result, HB formed the basis of a NNM-based model updating strategy and of convergence studies for nonlinear reduced order models \cite{kuether2014}. Nonlinear Fourier-based modal analysis was adapted to non-conservative autonomous dynamical systems in \cite{laxalde2009} thanks to generalized Fourier series with slow and fast time scales. Through this approach, refined later in \cite{krack2013}, NNMs of damped systems could be constructed and were found to give an accurate approximation of nonlinear resonances. By extending the spectral basis to two incommensurate frequencies, quasiperiodic (QP) oscillations can be studied, e.g., for monoharmonic \cite{schilder2005,peletan2014} and bi-periodic excitations \cite{guskov2012,liao2014}. Another interesting feature of nonlinear systems, namely the presence of co-existing periodic solutions and detached resonance curves, was studied in \cite{grolet2013} by coupling HB with Groebner bases.

The stability of nonlinear solutions can be studied by embedding Floquet theory within the HB formalism \cite{vongroll2001}, which is referred to as Hill's method \cite{hill1886}. In \cite{lanza2007}, Lanza \textit{et al.} provided an analytical approximation of Floquet exponents using the DF method. A semi-analytical version was developed by Bonani and Gilli for an arbitrary number of harmonics \cite{bonani1999}. However, these developments are limited to systems expressed in Lur'e form. Recently, more general theories have been proposed to compute Floquet exponents for autonomous \cite{traversa2008} and nonautonomous \cite{lazarus2010,traversa2012} systems.

Beside stability, bifurcations play a key role in the analysis of nonlinear systems. For example, fold bifurcations translate into a stability change of the periodic solutions, whereas QP oscillations are created or eliminated through Neimark-Sacker (NS) bifurcations. However, even if there exists a large body of literature on HB applied to nonlinear systems, very few studies attempted to use the method for tracking bifurcations.
In \cite{lanza2007}, bifurcation tracking is limited to single-harmonic approximations. Piccardi \cite{piccardi1994} proposed a procedure to obtain flip and fold curves (and even conditions for codimension-2 bifurcations \cite{piccardi1996}), but this procedure cannot describe NS bifurcations. Traversa \textit{et al.} developed a bifurcation tracking technique adapted to fold, flip and NS bifurcations \cite{traversa2008} by appending to the HB equation system an equation which describes the considered bifurcations through the Floquet multipliers. However, the resolution of the extra equation with the secant method makes the implementation inefficient when the size of the system increases.

In this context, the main contribution of the present paper is to adapt classical tools for bifurcation analysis in codimension-2 parameter space to the HB formalism. Because we target large-scale mechanical systems with localized nonlinearities, an algorithm that efficiently combines the computation of the Floquet exponents and bordering techniques is developed. A new procedure for the tracking of NS bifurcations that exploits the properties of eigenvalue derivatives is also proposed.

The paper is organized as follows. Section 2 recalls the theory of HB and its formulation in the framework of a continuation algorithm. In Section 3, Hill's method is introduced for assessing the stability of periodic solutions and for detecting their bifurcations. The proposed bifurcation tracking procedure is then presented with its adaptations for fold, branch point and NS bifurcations. The overall methodology is demonstrated using numerical experiments of a spacecraft structure that possesses a nonlinear vibration isolation device. Finally, the conclusions of the present study are summarized in Section 5.

%%%%%%%%%%%%%%%%%%%%%%%%%%%
% ---------------- SECTION 1
%%%%%%%%%%%%%%%%%%%%%%%%%%%
\section{Harmonic balance for periodic solutions}\label{HBmethod}

\subsection{Formulation of the dynamics in the frequency domain}\label{HB_theory}

We consider nonautonomous nonlinear dynamical systems with $n$ DOFs governed by the equations of motion
\medskip
\begin{equation}
\textbf{M}\ddot{\textbf{x}}+\textbf{C}\dot{\textbf{x}}+\textbf{K}\textbf{x} + \textbf{f}_{nl}(\textbf{x},\dot{\textbf{x}}) = \textbf{f}_{ext}(\omega,t)
\label{eom}
\end{equation}

where $\textbf{M}$, $\textbf{C}$ and $\textbf{K}$ are the mass, damping and stiffness matrices, respectively. Vectors $\textbf{x}$, $\textbf{f}_{nl}$ and $\textbf{f}_{ext}$ represent the displacements, the nonlinear forces and  the periodic external forces that are considered to be harmonic with frequency $\omega$ herein. The dots refer to the derivatives with respect to time $t$.

The periodic signals $\textbf{x}(t)$ and $\textbf{f}(\textbf{x},\dot{\textbf{x}},\omega,t)=\textbf{f}_{ext}(\omega,t)-\textbf{f}_{nl}(\textbf{x},\dot{\textbf{x}})$ in Equation (\ref{eom}) are approximated by Fourier series truncated to the $N_{H}$-th harmonic:
\medskip
\begin{eqnarray}
\textbf{x}(t) &=& \frac{\textbf{c}^{x}_{0}}{\sqrt{2}} + \sum_{k = 1}^{N_{H}}{\left(\textbf{s}^{x}_{k}\sin\left(\frac{k\omega t}{\nu}\right)+\textbf{c}^{x}_{k}\cos\left(\frac{k\omega t}{\nu}\right)\right)}\label{seriex}\\ \nonumber \\
\textbf{f}(t) &=& \frac{\textbf{c}^{f}_{0}}{\sqrt{2}} + \sum_{k = 1}^{N_{H}}{\left(\textbf{s}^{f}_{k}\sin\left(\frac{k\omega t}{\nu}\right)+\textbf{c}^{f}_{k}\cos\left(\frac{k\omega t}{\nu}\right)\right)}
\label{serie}
\end{eqnarray}

where $\textbf{s}_{k}$ and $\textbf{c}_{k}$ represent the vectors of the Fourier coefficients related to the sine and cosine terms, respectively, and the integer $\nu$ accounts for subharmonics of the excitation frequency $\omega$. The Fourier coefficients of $\textbf{f}\left(t\right)$, $\textbf{c}^{f}_{k}$ and $\textbf{s}^{f}_{k}$, depend on the Fourier coefficients of the displacements $\textbf{x}\left(t\right)$, $\textbf{c}^{x}_{k}$ and $\textbf{s}^{x}_{k}$, which represent the new unknowns of the problem. These coefficients are gathered into the $\left(2N_{H}+1\right)n \times 1$ vectors
\medskip
\begin{equation}
\textbf{z} = \left[\begin{array}{cccccc}\left(\textbf{c}^{x}_{0}\right)^T & \left(\textbf{s}^{x}_{1}\right)^T & \left(\textbf{c}^{x}_{1}\right)^T & \hdots & \left(\textbf{s}^{x}_{N_{H}}\right)^T & \left(\textbf{c}^{x}_{N_{H}}\right)^T \end{array}\right]^T
\end{equation}
\begin{equation}
\textbf{b} = \left[\begin{array}{cccccc}\left(\textbf{c}^{f}_{0}\right)^T & \left(\textbf{s}^{f}_{1}\right)^T & \left(\textbf{c}^{f}_{1}\right)^T & \hdots & \left(\textbf{s}^{f}_{N_{H}}\right)^T & \left(\textbf{c}^{f}_{N_{H}}\right)^T \end{array}\right]^T
\end{equation}

The displacements and forces are recast into a more compact form
\medskip
\begin{eqnarray}
\textbf{x}(t) &=& \left(\textbf{Q}(t) \otimes \mathbb{I}_n\right)\textbf{z}\label{approx_serie_dep}\\
\textbf{f}(t) &=& \left(\textbf{Q}(t) \otimes \mathbb{I}_n\right)\textbf{b}\label{approx_serie_f}
\end{eqnarray}

where $\otimes$ and $\mathbb{I}_n$ stand for the Kronecker tensor product and the identity matrix of size $n$, respectively, and $\textbf{Q}(t)$ is a matrix containing the sine and cosine series
\medskip
\begin{equation}
\textbf{Q}(t) = \left[\begin{array}{cccccc}\frac{1}{\sqrt{2}} & \sin\left(\frac{\omega t}{\nu}\right) & \cos\left(\frac{\omega t}{\nu}\right) & \hdots & \sin\left(N_H\frac{\omega t}{\nu}\right) & \cos\left(N_H\frac{\omega t}{\nu}\right) \end{array}\right]
\end{equation}

Velocities and accelerations can also be defined using the Fourier series \cite{jaumouille2010}, with
\medskip
\begin{eqnarray}
\dot{\textbf{x}}(t) &=& \left(\dot{\textbf{Q}}(t) \otimes \mathbb{I}_n\right)\textbf{z}=\left(\left(\textbf{Q}(t)\boldsymbol{\nabla}\right) \otimes \mathbb{I}_n\right)\textbf{z}\label{approx_serie_vel}\\
\ddot{\textbf{x}}(t) &=& \left(\ddot{\textbf{Q}}(t) \otimes \mathbb{I}_n\right)\textbf{z}=\left(\left(\textbf{Q}(t)\boldsymbol{\nabla}^2\right) \otimes \mathbb{I}_n\right)\textbf{z}\label{approx_serie_acc}
\end{eqnarray}
where
\medskip
\begin{equation}
\begin{array}{ccc}
\boldsymbol{\nabla} = \left[\begin{array}{ccccc}0 & & & & \\ & \ddots & & & \\ & & \boldsymbol{\nabla}_k & & \\ & & & \ddots & \\ & & & & \boldsymbol{\nabla}_{N_H}\end{array}\right], & & \boldsymbol{\nabla}\boldsymbol{\nabla} = \boldsymbol{\nabla}^2 = \left[\begin{array}{ccccc}0 & & & & \\ & \ddots & & & \\ & & \boldsymbol{\nabla}^2_k & & \\ & & & \ddots & \\ & & & & \boldsymbol{\nabla}^2_{N_H}\end{array}\right]
\end{array}
\end{equation}

with
\medskip
\begin{equation}
\begin{array}{ccc}
\boldsymbol{\nabla}_k = \left[\begin{array}{cc}0 & -\frac{k\omega}{\nu}\\ \frac{k\omega}{\nu} & 0\end{array}\right], & & \boldsymbol{\nabla}^2_k = \left[\begin{array}{cc}-\left(\frac{k\omega}{\nu}\right)^2 & 0\\ 0 & -\left(\frac{k\omega}{\nu}\right)^2\end{array}\right]
\end{array}
\end{equation}

Substituting expressions (\ref{approx_serie_dep})-(\ref{approx_serie_f}) and (\ref{approx_serie_vel})-(\ref{approx_serie_acc}) in the equations of motion (\ref{eom}) yields
\medskip
\begin{equation}
\textbf{M}\left(\left(\textbf{Q}(t)\boldsymbol{\nabla}^2\right) \otimes \mathbb{I}_n\right)\textbf{z} + \textbf{C}\left(\left(\textbf{Q}(t)\boldsymbol{\nabla}\right) \otimes \mathbb{I}_n\right)\textbf{z}+\textbf{K}\left(\textbf{Q}(t) \otimes \mathbb{I}_n\right)\textbf{z} = \left(\textbf{Q}(t) \otimes \mathbb{I}_n\right)\textbf{b} \label{eom_v2}
\end{equation}

Considering the mixed-product property of the Kronecker tensor product $\left(\textbf{A} \otimes \textbf{B}\right)\left(\textbf{C} \otimes \textbf{D}\right) = \left(\textbf{A}\textbf{C}\right) \otimes \left(\textbf{B}\textbf{D}\right)$ yields
\medskip
\begin{eqnarray}
\textbf{M}\left(\left(\textbf{Q}(t)\boldsymbol{\nabla}^2\right) \otimes \mathbb{I}_n\right) &=& \left(1 \otimes \textbf{M}\right)\left(\left(\textbf{Q}(t)\boldsymbol{\nabla}^2\right) \otimes \mathbb{I}_n\right) = \left(\textbf{Q}(t)\boldsymbol{\nabla}^2\right) \otimes \textbf{M} \\
\textbf{C}\left(\left(\textbf{Q}(t)\boldsymbol{\nabla}\right) \otimes \mathbb{I}_n\right) &=& \left(1 \otimes \textbf{C}\right)\left(\left(\textbf{Q}(t)\boldsymbol{\nabla}\right) \otimes \mathbb{I}_n\right) = \left(\textbf{Q}(t)\boldsymbol{\nabla}\right) \otimes \textbf{C} \\
\textbf{K}\left(\textbf{Q}(t) \otimes \mathbb{I}_n\right) &=& \left(1 \otimes \textbf{K}\right)\left(\textbf{Q}(t) \otimes \mathbb{I}_n\right) = \textbf{Q}(t)\otimes \textbf{K}
\end{eqnarray}

These expressions are plugged into Equation (\ref{eom_v2}), which gives
\medskip
\begin{equation}
\left(\left(\textbf{Q}(t)\boldsymbol{\nabla}^2\right) \otimes \textbf{M}\right)\textbf{z}+\left(\left(\textbf{Q}(t)\boldsymbol{\nabla}\right) \otimes \textbf{C}\right)\textbf{z}+\left(\textbf{Q}(t)\otimes \textbf{K}\right)\textbf{z} = \left(\textbf{Q}(t) \otimes \mathbb{I}_n\right)\textbf{b} \label{eom_v3}
\end{equation}

In order to remove the time dependency and to obtain an expression relating the different Fourier coefficients, a Galerkin procedure projects Equation (\ref{eom_v3}) on the orthogonal trigonometric basis $\textbf{Q}(t)$
\medskip
\begin{multline}
\left(\left(\frac{2}{T}\int_{0}^{T}{\textbf{Q}^T(t)\textbf{Q}(t)\,dt}\,\boldsymbol{\nabla}^2\right) \otimes \textbf{M}\right)\textbf{z}+\left(\left(\frac{2}{T}\int_{0}^{T}{\textbf{Q}^T(t)\textbf{Q}(t)\,dt}\,\boldsymbol{\nabla}\right) \otimes \textbf{C}\right)\textbf{z}+\hdots \\ \left(\left(\frac{2}{T}\int_{0}^{T}{\textbf{Q}^T(t)\textbf{Q}(t)\,dt}\right)\otimes \textbf{K}\right)\textbf{z} = \left(\left(\frac{2}{T}\int_{0}^{T}{\textbf{Q}^T(t)\textbf{Q}(t)\,dt}\right) \otimes \mathbb{I}_n\right)\textbf{b}\label{eom_v4}
\end{multline}

where $T$ is the period of the external force. Considering that
\medskip
\begin{equation}
\frac{2}{T}\int_{0}^{T}{\textbf{Q}^T(t)\textbf{Q}(t)dt} = \mathbb{I}_{2N_H+1}
\end{equation}

the equations of motion expressed in the frequency domain are eventually obtained
\begin{equation}
\left(\boldsymbol{\nabla}^2\otimes \textbf{M}\right)\textbf{z}+\left(\boldsymbol{\nabla}\otimes \textbf{C}\right)\textbf{z}+\left(\mathbb{I}_{2N_H+1}\otimes \textbf{K}\right)\textbf{z} = \left(\mathbb{I}_{2N_H+1} \otimes \mathbb{I}_n\right)\textbf{b}
\end{equation}

or, in a more compact form,
\begin{equation}
\textbf{h}(\textbf{z},\omega) \equiv \textbf{A}(\omega)\textbf{z}-\textbf{b}(\textbf{z}) = \textbf{0}
\label{eomf}
\end{equation}

where $\textbf{A}$ is the $\left(2N_{H}+1\right)n \times \left(2N_{H}+1\right)n$ matrix describing the linear dynamics
\medskip
\begin{eqnarray}
\textbf{A} &=& \nonumber \boldsymbol{\nabla}^2\otimes \textbf{M}+\boldsymbol{\nabla}\otimes \textbf{C}+\mathbb{I}_{2N_H+1}\otimes \textbf{K}\\
&=& \left[ \begin{array}{cccc}
\textbf{K} & & & \\
 & \begin{array}{cc}\textbf{K}-\left(\frac{\omega}{\nu}\right)^{2}\textbf{M} & -\frac{\omega}{\nu} \textbf{C}\\ \frac{\omega}{\nu} \textbf{C} & \textbf{K}-\left(\frac{\omega}{\nu}\right)^{2}\textbf{M}\end{array} & & \\
 & & \ddots & \\
 & & & \begin{array}{cc}\textbf{K}-\left(N_{H}\frac{\omega}{\nu}\right)^{2}\textbf{M} & -N_{H}\frac{\omega}{\nu} \textbf{C}\\ N_{H}\frac{\omega}{\nu} \textbf{C} & \textbf{K}-\left(N_{H}\frac{\omega}{\nu}\right)^{2}\textbf{M}\end{array} \\
\end{array}
\right] \nonumber \\
& &
\end{eqnarray}

Expression (\ref{eomf}) can be seen as the equations of amplitude of (\ref{eom}), i.e., if $\textbf{z}^{*}$ is a root of (\ref{eomf}), then the time signals $\textbf{x}^{*}$ constructed from $\textbf{z}^{*}$ with (\ref{seriex}) are solutions of (\ref{eom}).

\subsection{Analytical expression of the nonlinear terms and of the jacobian matrix of the system}

Equation (\ref{eomf}) is nonlinear and has to be solved iteratively (e.g., with a Newton-Raphson procedure, or with the hybrid Powell nonlinear solver \cite{powell1970,nacivet2003}). At each iteration, an evaluation of $\textbf{b}$ and of $\partial \textbf{h}/\partial \textbf{z}$ has to be provided. When $\textbf{f}$ can be accurately approximated with a few number of harmonics and when its analytical sinusoidal expansion is known \cite{lee2005}, or for some types of restoring force \cite{petrov2003}, analytical expressions relating the Fourier coefficients of the forces $\textbf{b}$ and of the displacements $\textbf{z}$ can be obtained together with the expression of the Jacobian matrix of the system. Otherwise, the alternating frequency/time-domain (AFT) technique \cite{cameron1989} can be used to compute $\textbf{b}$:
\medskip
\begin{equation}
\textbf{z} \xrightarrow[]{\mbox{FFT}^{-1}} \textbf{x}(t) \xrightarrow[]{} \textbf{f}\left(\textbf{x},\dot{\textbf{x}},\omega,t\right) \xrightarrow[]{\mbox{FFT}} \textbf{b}(\textbf{z})
\end{equation}

The Jacobian matrix of the system can then be computed through finite differences, which is computationally demanding.

An efficient alternative consists in rewriting the inverse Fourier transform as a linear operator $\boldsymbol{\Gamma}\left(\omega\right)$. First proposed by Hwang \cite{Hwang1991}, the method, also called \textit{trigonometric collocation}, was adapted to the AFT technique \cite{Xie1996}, and has been widely used since then \cite{bonani1999,narayanan1998,duan2005,kim2005}. Denoting $N$ as the number of time samples of a discretized period of oscillation, one defines vectors $\tilde{\textbf{x}}$ and $\tilde{\textbf{f}}$ containing the concatenated $nN$ time samples of the displacements and the forces, respectively, for all DOFs:
\medskip
\begin{eqnarray}
\tilde{\textbf{x}} &=& \left[\begin{array}{ccccccc}x_{1}\left(t_{1}\right) & \hdots & x_{1}\left(t_{N}\right) & \hdots & x_{n}\left(t_{1}\right) & \hdots & x_{n}\left(t_{N}\right) \end{array}\right]^T\\
\tilde{\textbf{f}} &=& \left[\begin{array}{ccccccc}f_{1}\left(t_{1}\right) & \hdots & f_{1}\left(t_{N}\right) & \hdots & f_{n}\left(t_{1}\right) & \hdots & f_{n}\left(t_{N}\right) \end{array}\right]^T
\end{eqnarray}

The inverse Fourier transform can then be written as a linear operation:
\medskip
\begin{equation}
\tilde{\textbf{x}} = \boldsymbol{\Gamma}\left(\omega\right)\textbf{z}
\end{equation}

with the $nN \times \left(2N_{H}+1\right)n$ sparse operator
\medskip
\begin{multline}
\boldsymbol{\Gamma}\left(\omega\right) = \left[\begin{array}{cccc}\mathbb{I}_{n}\otimes \left[\begin{array}{c} 1/\sqrt{2} \\ 1/\sqrt{2} \\ \vdots \\ 1/\sqrt{2} \end{array} \right] & \mathbb{I}_{n}\otimes \left[\begin{array}{c} \sin\left(\frac{\omega t_{1}}{\nu}\right) \\ \sin\left(\frac{\omega t_{2}}{\nu}\right) \\ \vdots \\ \sin\left(\frac{\omega t_{N}}{\nu}\right) \end{array} \right] & \mathbb{I}_{n}\otimes \left[\begin{array}{c} \cos\left(\frac{\omega t_{1}}{\nu}\right) \\ \cos\left(\frac{\omega t_{2}}{\nu}\right) \\ \vdots \\ \cos\left(\frac{\omega t_{N}}{\nu}\right) \end{array} \right] & \hdots \end{array}\right. \\ \\ \left. \begin{array}{cc} \mathbb{I}_{n}\otimes \left[\begin{array}{c} \sin\left(N_{H}\frac{\omega t_{1}}{\nu}\right) \\ \sin\left(N_{H}\frac{\omega t_{2}}{\nu}\right) \\ \vdots \\ \sin\left(N_{H}\frac{\omega t_{N}}{\nu}\right) \end{array} \right] & \mathbb{I}_{n}\otimes \left[\begin{array}{c} \cos\left(N_{H}\frac{\omega t_{1}}{\nu}\right) \\ \cos\left(N_{H}\frac{\omega t_{2}}{\nu}\right) \\ \vdots \\ \cos\left(N_{H}\frac{\omega t_{N}}{\nu}\right) \end{array} \right] \end{array} \right]
\end{multline}

\begin{figure}[h!t]
\centering \includegraphics[scale=0.8]{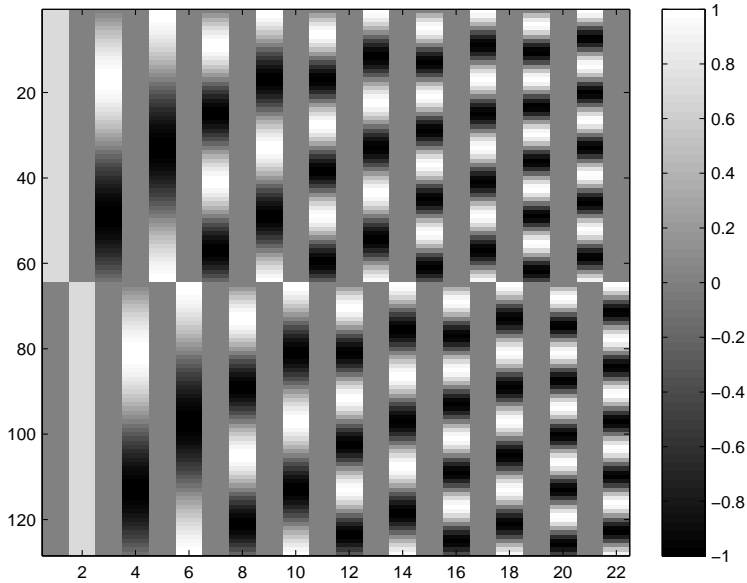}
\caption{Illustration of the inverse Fourier transformation matrix $\boldsymbol{\Gamma}\left(\omega\right)$ for $n = 2$, $N = 64$, and $N_{H} = 5$.} \label{fourier_mat}
\end{figure}

Figure \ref{fourier_mat} represents the inverse transformation matrix for the case $n = 2$, $N = 64$, and $N_{H} = 5$.

The direct Fourier transformation is written
\medskip
\begin{equation}
\textbf{z} = \left(\boldsymbol{\Gamma}\left(\omega\right)\right)^{+}\tilde{\textbf{x}}
\end{equation}

where $^{+}$ stands for the Moore-Penrose pseudoinverse $\boldsymbol{\Gamma}^{+} = \boldsymbol{\Gamma}^{T}\left(\boldsymbol{\Gamma}\boldsymbol{\Gamma}^{T}\right)^{-1}$. The Fourier coefficients of the external and nonlinear forces are simply obtained by transforming the signals in the time domain back to the frequency domain:
\medskip
\begin{equation}
\textbf{b}(\textbf{z}) = \left(\boldsymbol{\Gamma}\left(\omega\right)\right)^{+}\tilde{\textbf{f}}
\end{equation}

The Jacobian matrix of expression (\ref{eomf}) with respect to the Fourier coefficients $\textbf{z}$, denoted as $\textbf{h}_{z}$, can be computed as in \cite{bonani1999,narayanan1998,kim2005}
\medskip
\begin{equation}
\textbf{h}_{z} = \frac{\partial\textbf{h}}{\partial\textbf{z}} = \textbf{A} - \frac{\partial\textbf{b}}{\partial\textbf{z}} = \textbf{A} - \frac{\partial\textbf{b}}{\partial\tilde{\textbf{f}}}\frac{\partial\tilde{\textbf{f}}}{\partial\tilde{\textbf{x}}}\frac{\partial\tilde{\textbf{x}}}{\partial\textbf{z}}=\textbf{A} - \boldsymbol{\Gamma}^{+}\frac{\partial\tilde{\textbf{f}}}{\partial\tilde{\textbf{x}}}\boldsymbol{\Gamma}
\end{equation}

In general, the derivatives of the forces with respect to the displacements in the time domain can be expressed analytically, which leads to a very effective computation of the Jacobian matrix. These derivatives have to be evaluated only for the nonlinear DOFs.

\subsection{Continuation procedure}

Usually, the frequency response is to be computed in a range of frequencies, rather than for a single frequency $\omega$. In this paper, a continuation procedure based on tangent predictions and Moore-Penrose corrections, as in the software \textsc{matcont} \cite{dhooge2003}, is coupled to HB. The search for a tangent vector $\textbf{t}^{(i)}$ at an iteration point $\left(\textbf{z}^{(i)},\omega^{(i)}\right)$ along the branch reads
\medskip
\begin{equation}
\left[\begin{array}{cc}\textbf{h}_{z} & \textbf{h}_{\omega} \\  \multicolumn{2}{c}{\textbf{t}^{T}_{(i-1)}}\end{array}\right] \textbf{t}_{(i)} = \left[\begin{array}{c}\textbf{0} \\ 1\end{array}\right]
\label{tangent}
\end{equation}

where $\textbf{h}_{\omega}$ stands for the derivative of $\textbf{h}$ with respect to $\omega$,
\medskip
\begin{equation}
\textbf{h}_{\omega} = \frac{\partial \textbf{h}}{\partial \omega} = \frac{\partial \textbf{A}}{\partial \omega}\textbf{z}
\end{equation}

The last equation in (\ref{tangent}) prevents the continuation procedure from turning back. For the first iteration of the procedure, the sum of the components of the tangent is imposed to be equal to 1.

The correction stage is based on Newton's method. Introducing new optimization variables $\textbf{v}_{(i,j)}$ initialized as $\textbf{v}_{(i,1)} = \textbf{t}_{(i)}$, and $\textbf{y}_{(i,j)} = \left[\textbf{z}_{(i,j)} \ \omega_{(i,j)}\right]^{T}$, the Newton iterations are constructed as follows:
\medskip
\begin{equation}
\begin{array}{l}
\textbf{y}_{(i,j+1)} = \textbf{y}_{(i,j)} - \textbf{G}_{y}^{-1}\left(\textbf{y}_{(i,j)},\textbf{v}_{(i,j)}\right)\textbf{G}\left(\textbf{y}_{(i,j)},\textbf{v}_{(i,j)}\right)\\
\textbf{v}_{(i,j+1)} = \textbf{v}_{(i,j)} - \textbf{G}_{y}^{-1}\left(\textbf{y}_{(i,j)},\textbf{v}_{(i,j)}\right)\textbf{R}\left(\textbf{y}_{(i,j)},\textbf{v}_{(i,j)}\right)\\
\end{array}
\label{corr_cont1}
\end{equation}

with
\medskip
\begin{equation}
\begin{array}{c}
\begin{array}{ccc}\textbf{G}\left(\textbf{y},\textbf{v}\right) = \left[\begin{array}{c} \textbf{h}\left(\textbf{y}\right)\\ \textbf{0}\end{array}\right], & & \textbf{G}_{y}\left(\textbf{y},\textbf{v}\right) = \left[\begin{array}{c} \begin{array}{cc}\textbf{h}_{z}\left(\textbf{y}\right) & \textbf{h}_{\omega}\left(\textbf{y}\right)\end{array}\\ \textbf{v}^{T}\end{array}\right],\end{array} \\ \\ \textbf{R}\left(\textbf{y},\textbf{v}\right) = \left[\begin{array}{c} \left[\begin{array}{cc}\textbf{h}_{z}\left(\textbf{y}\right) & \textbf{h}_{\omega}\left(\textbf{y}\right)\end{array}\right]\textbf{v}\\ \textbf{0}\end{array}\right]
\end{array}
\label{corr_cont2}
\end{equation}

%%%%%%%%%%%%%%%%%%%%%%%%%%%%%
%% ---------------- SECTION 2
%%%%%%%%%%%%%%%%%%%%%%%%%%%%%
\section{Harmonic balance for bifurcation detection and tracking}

\subsection{Stability analysis}\label{stab_anal_sec}

The continuation procedure presented in Section \ref{HBmethod} does not indicate if a periodic solution is stable or not. In the case of time-domain methods, such as the shooting technique, a by-product of the continuation procedure is the monodromy matrix of the system \cite{peeters2009}. Its eigenvalues are termed the \textit{Floquet multipliers} from which the stability of the solution can be deduced. For frequency-domain techniques, one can use a variant of the Floquet theory, the so-called \textit{Hill's method}, whose coefficients are also obtained as a by-product of the calculations. Hill's method is known to give accurate results for a reasonable number of harmonics $N_H$, and to be effective for large systems \cite{peletan2013}.

Following the procedure described in \cite{vongroll2001}, a periodic solution $\textbf{x}^*(t)$ satisfying (\ref{eom}) is perturbed with a periodic solution $\textbf{s}(t)$ modulated by an exponential decay:
\medskip
\begin{equation}
\textbf{p}(t) = \textbf{x}^*(t)+e^{\lambda t}\textbf{s}(t)
\label{perturb}
\end{equation}

Introducing this perturbation into expression (\ref{eom}) yields
\medskip
\begin{equation}
\textbf{M}\ddot{\textbf{x}}^*+\textbf{C}\dot{\textbf{x}}^*+\textbf{K}\textbf{x}^*+\left(\lambda^2\textbf{M}\textbf{s}+\lambda\left(2\textbf{M}\dot{\textbf{s}}+\textbf{C}\textbf{s}\right)+\textbf{M}\ddot{\textbf{s}}+\textbf{C}\dot{\textbf{s}}+\textbf{K}\textbf{s}\right)e^{\lambda t}=\textbf{f}\left(\textbf{p},\dot{\textbf{p}},\omega,t\right)
\end{equation}

By approximating the solution and the perturbation as Fourier series truncated to the $N_H$-th order, i.e., $\textbf{x}^*(t)= \left(\textbf{Q}(t) \otimes \mathbb{I}_n\right)\textbf{z}^*$  and $\textbf{s}(t)=\left(\textbf{Q}(t) \otimes \mathbb{I}_n\right)\textbf{u}$, and by applying a Galerkin procedure as in Section \ref{HB_theory}, we obtain
\medskip
\begin{equation}
\textbf{A}\textbf{z}^*+\left(\boldsymbol{\Delta}_{2}\lambda^2+\boldsymbol{\Delta}_{1}\lambda+\textbf{A}\right)e^{\lambda t}\textbf{u} = \textbf{b}\left(\textbf{z}^*+e^{\lambda t}\textbf{u}\right)
\label{hill_prev}
\end{equation}

with
\begin{eqnarray}
\boldsymbol{\Delta}_{1} &=& \boldsymbol{\nabla}\otimes 2\textbf{M}+\mathbb{I}_{2N_H+1}\otimes \textbf{C} \\
&=& \left[\begin{array}{cccc}
\textbf{C} & & & \\
 & \begin{array}{cc}\textbf{C} & -2\frac{\omega}{\nu}\textbf{M}\\ 2\frac{\omega}{\nu}\textbf{M} & \textbf{C}\end{array} & & \\
 & & \ddots & \\
 & & & \begin{array}{cc}\textbf{C} & -2N_{H}\frac{\omega}{\nu}\textbf{M}\\ 2N_{H}\frac{\omega}{\nu}\textbf{M} & \textbf{C}\end{array}
\end{array}
\right]\nonumber\\
\boldsymbol{\Delta}_{2} &=& \mathbb{I}_{2N_{H}+1} \otimes \textbf{M}
\end{eqnarray}

The right-hand side of Equation (\ref{hill_prev}) is evaluated through a Taylor series expansion around the solution $\textbf{z}^*$
\medskip
\begin{equation}
\textbf{b}\left(\textbf{z}^*+e^{\lambda t}\textbf{u}\right) = \textbf{b}\left(\textbf{z}^*\right)+\left.\frac{\partial\textbf{b}}{\partial\textbf{z}}\right|_{\textbf{z} = \textbf{z}^*}\left(e^{\lambda t}\textbf{u}\right)
\label{taylor_b}
\end{equation}

Since $\textbf{A}\textbf{z}^*-\textbf{b}\left(\textbf{z}^*\right)$ is equal to zero by definition, and given that
\medskip
\begin{equation}
\textbf{A}-\left.\frac{\partial\textbf{b}}{\partial\textbf{z}}\right|_{\textbf{z} = \textbf{z}^*} = \textbf{h}_z
\end{equation}

replacing (\ref{taylor_b}) in (\ref{hill_prev}) yields
\medskip
\begin{equation}
\left(\boldsymbol{\Delta}_{2}\lambda^2+\boldsymbol{\Delta}_{1}\lambda+\textbf{h}_z\right)e^{\lambda t}\textbf{u} = \textbf{0}
\end{equation}

Hill's coefficients $\boldsymbol{\lambda}$ are thus the solutions of the quadratic eigenvalue problem
\medskip
\begin{equation}
\boldsymbol{\Delta}_{2}\lambda^{2}+\boldsymbol{\Delta}_{1}\lambda+\textbf{h}_{z} = \textbf{0}
\label{quad_pvp}
\end{equation}

When embedded in a continuation scheme, $\textbf{h}_{z}$ is already obtained from the correction stage (\ref{corr_cont1}-\ref{corr_cont2}). Since $\boldsymbol{\Delta}_{1}$ and $\boldsymbol{\Delta}_{2}$ are easily computed, the main computational effort amounts to solving a quadratic eigenvalue problem, which can be rewritten as a linear eigenvalue problem of double size
\medskip
\begin{equation}
\textbf{B}_{1}-\gamma\textbf{B}_{2} = \textbf{0}\label{eigenprob}
\end{equation}
with
\begin{equation}
\begin{array}{ccc}
\textbf{B}_{1} = \left[\begin{array}{cc}\boldsymbol{\Delta}_{1} & \textbf{h}_{z}\\ -\mathbb{I} & \textbf{0}\end{array}\right], & & \textbf{B}_{2} = -\left[\begin{array}{cc}\boldsymbol{\Delta}_{2} & \textbf{0} \\ \textbf{0} & \mathbb{I}\end{array}\right]
\end{array}\label{eqB1}
\end{equation}

The coefficients $\boldsymbol{\lambda}$ are found among the eigenvalues of the $\left(2N_{H}+1\right)2n \times \left(2N_{H}+1\right)2n$ matrix
\medskip
\begin{eqnarray}
\textbf{B} &=& \textbf{B}_{2}^{-1}\textbf{B}_{1}\label{eqB2}\\\label{DerivB}
&=& \left[\begin{array}{cc}-\boldsymbol{\Delta}_{2}^{-1}\boldsymbol{\Delta}_{1} & -\boldsymbol{\Delta}_{2}^{-1}\textbf{h}_{z}\\ \mathbb{I} & \textbf{0}\end{array}\right]
\end{eqnarray}

However, only $2n$ eigenvalues among the complete set $\boldsymbol{\lambda}$ approximate the Floquet exponents $\tilde{\boldsymbol{\lambda}}$ of the solution $\textbf{x}^*$ \cite{lazarus2010}. The other eigenvalues are spurious and do not have any physical meaning; their number also increases with the number of harmonics $N_H$. Moore \cite{moore2005} showed that the Floquet exponents $\tilde{\boldsymbol{\lambda}}$ are the $2n$ eigenvalues with the smallest imaginary part in modulus. The diagonal matrix
\medskip
\begin{equation}
\tilde{\textbf{B}} = \left[
\begin{array}{cccc}
\tilde{\lambda}_1 & & & \\
& \tilde{\lambda}_2 & & \\
& & \ddots & \\
& & & \tilde{\lambda}_{2n}
\end{array}\right]
\end{equation}
gathers the Floquet exponents and will play a key role for the detection and tracking of bifurcations in Sections \ref{bifdetection} and \ref{biftracking}.

Eventually, the stability of a periodic solution can be assessed, i.e., if at least one of the Floquet exponents has a positive real part, then the solution is unstable, otherwise it is asymptotically stable.

\subsection{Detection of bifurcations}\label{bifdetection}

To detect bifurcations, \textit{test functions} $\phi$ are evaluated at each iteration of the numerical continuation process \cite{dhooge2003}. The roots of these test functions indicate the presence of bifurcations.

A fold bifurcation is simply detected when the $i_{\omega}$-th component of the tangent prediction related to the active parameter $\omega$ changes sign \cite{govaerts2011}. A suitable test function is thus
\medskip
\begin{equation}
\phi_{F} = \textbf{t}_{i_{\omega}}
\end{equation}

According to its algebraic definition \cite{seydel2010}, we note that a fold bifurcation is characterized by a rank deficiency of 1 of the Jacobian matrix $\textbf{h}_{z}$, with $\textbf{h}_{\omega} \notin \text{range}\left(\textbf{h}_{z}\right)$. Another (more computationally demanding) test function is therefore
\medskip
\begin{equation}\label{Fold2}
\phi_{F} = \left|\textbf{h}_z\right|
\end{equation}

Similarly, the Jacobian matrix is rank deficient for branch point (BP) bifurcations, with $\textbf{h}_{\omega} \in \text{range}\left(\textbf{h}_{z}\right)$. They can be detected using the test function for fold bifurcations (\ref{Fold2}), but a more specific test function is \cite{govaerts2011}
\medskip
\begin{equation}\label{BPtest}
\begin{array}{l}
\phi_{BP} = \left|\begin{array}{cc}\textbf{h}_{z} & \textbf{h}_{\omega} \\  \multicolumn{2}{c}{\textbf{t}^{T}}\end{array}\right|
\end{array}
\end{equation}

The Neimark-Sacker (NS) bifurcation is the third bifurcation studied in this paper. It is detected when a pair of Floquet exponents crosses the imaginary axis as a pair of complex conjugates. According to \cite{guckenheimer1997}, the bialternate product of a $m \times m$ matrix $\textbf{P}$, $\textbf{P}_{\odot} = \textbf{P} \odot \mathbb{I}_{m}$, has a dimension $m(m-1)/2$ and has the property to be singular when two eigenvalues of $\textbf{P}$ are two purely imaginary complex conjugates. As a result, the test function for NS bifurcations is
\medskip
\begin{equation}
\phi_{NS} = \left|\tilde{\textbf{B}}_{\odot}\right|
\end{equation}

The size of the bialternate product rapidly increases with the number of DOFs $n$, but the diagonal shape of $\tilde{\textbf{B}}$ implies that $\tilde{\textbf{B}}_{\odot}$ is also diagonal, which allows for a fast evaluation of its terms.

To overcome the issue of computing determinants for large-scale systems, the so-called \textit{bordering technique} replaces the evaluation of the determinant of a matrix $\textbf{G}$ with the evaluation of a scalar function $g$ which vanishes at the same time as the determinant \cite{beyn2002}. The function $g$ is obtained by solving the bordered system
\medskip
\begin{equation}\label{BorderedS}
\left[\begin{array}{cc}\textbf{G} & \textbf{p}\\\textbf{q}^{*} & 0\end{array}\right]\left[\begin{array}{c}\textbf{w}\\g\end{array}\right] = \left[\begin{array}{c}\textbf{0}\\1\end{array}\right]
\end{equation}

where $^{*}$ denotes conjugate transpose, and vectors $\textbf{p}$ and $\textbf{q}$ are chosen to ensure the nonsingularity of the system of equations. When $\textbf{G}$ is almost singular, $\textbf{p}$ and $\textbf{q}$ are chosen close to the null vectors of $\textbf{G}^{*}$ and $\textbf{G}$, respectively. For instance, for NS bifurcations $\textbf{G} = \tilde{\textbf{B}}_{\odot}$.

\subsection{Tracking of bifurcations}\label{biftracking}

Once a bifurcation is detected, it can be tracked with respect to an additional parameter. To continue codimension-1 bifurcations with respect to two parameters, such as, e.g., the frequency and amplitude of the external forcing, the equation defining the bifurcation $g=0$ is appended to (\ref{eomf})
\medskip
\begin{equation}\left\{
\begin{array}{l}
\textbf{h} \equiv \textbf{A}\textbf{z}-\textbf{b} = 0\vspace{2mm}\\
g = 0
\end{array}
\right.
\end{equation}

For fold and BP bifurcations, $g$ is the solution of the bordered system (\ref{BorderedS}) with $\textbf{G}=\textbf{h}_{z}$. For NS bifurcations, $\textbf{G}=\tilde{\textbf{B}}_{\odot}$ is considered in the bordered system.

During the continuation procedure, the computation of the derivatives of the additional equation is required. As shown in \cite{beyn2002}, analytical expressions for the derivatives of $g$ with respect to $\alpha$, where $\alpha$ denotes a component of $\textbf{z}$ or one of the two active parameters, are found as
\medskip
\begin{equation}
g_{\alpha} = -\textbf{v}^{*}\textbf{G}_{\alpha}\textbf{w}
\end{equation}

where $\textbf{G}_{\alpha}$ is the derivative of \textbf{G} with respect to $\alpha$, and where $\textbf{w}$ and $\textbf{v}$ comes from the resolution of the bordered system and its transposed version:
\medskip
\begin{equation}
\left[\begin{array}{cc}\textbf{G} & \textbf{p}\\\textbf{q}^{*} & 0\end{array}\right]\left[\begin{array}{c}\textbf{w}\\g\end{array}\right] = \left[\begin{array}{c}\textbf{0}\\1\end{array}\right]
\end{equation}

\begin{equation}
\left[\begin{array}{cc}\textbf{G} & \textbf{p}\\\textbf{q}^{*} & 0\end{array}\right]^{*}\left[\begin{array}{c}\textbf{v}\\e\end{array}\right] = \left[\begin{array}{c}\textbf{0}\\1\end{array}\right]
\end{equation}

As a result, the only term that has to be evaluated is $\textbf{G}_{\alpha}$. For fold and BP bifurcations,
\medskip
\begin{equation}
\textbf{G}_{\alpha} = \textbf{h}_{z\alpha}
\end{equation}

where $\textbf{h}_{z\alpha}$ is the derivative of the Jacobian $\textbf{h}_{z}$ with respect to $\alpha$ and is computed through finite differences.

For NS bifurcations,
\medskip
\begin{equation}\label{eigen_deriv}
\begin{array}{lll}
 \textbf{G}_{\alpha} & = & \frac{\partial}{\partial \alpha}\left(\tilde{\textbf{B}}_{\odot}\right) = \left(\frac{\partial}{\partial \alpha}\left(\tilde{\textbf{B}}\right)\right)_{\odot}=\left[
\begin{array}{cccc}
\frac{\partial\tilde{\lambda}_1}{\partial\alpha} & & & \\
& \frac{\partial\tilde{\lambda}_2}{\partial\alpha} & & \\
& & \ddots & \\
& & & \frac{\partial\tilde{\lambda}_{2n}}{\partial\alpha}
\end{array}\right]_{\odot}
\end{array}
\end{equation}

Finite differences could be used to compute the derivatives of the Floquet exponents. However, this means that the eigenvalue problem (\ref{eigenprob}) has to be solved for each perturbation of the components of $\textbf{z}$, and for the perturbation of the two continuation parameters. This represents a total of $n\left(2N_H+1\right)+2$ resolutions of the eigenvalue problem per iteration, which is cumbersome for large systems. Instead, we propose to compute the derivatives in (\ref{eigen_deriv}) using the properties of eigenvalue derivatives demonstrated by Van der Aa \textit{et al.} \cite{vanderaa2007}. Denoting as $\boldsymbol{\Lambda}$ the eigenvector matrix of $\textbf{B}$, and $\boldsymbol{\xi}$ the localization vector containing the index of the $2n$ Floquet exponents $\tilde{\boldsymbol{\lambda}}$ among the eigenvalues $\boldsymbol{\lambda}$, i.e., $\tilde{\lambda}_i = \lambda_{\xi_i}$, the eigenvalues derivatives can be computed as
\medskip
\begin{equation}
\frac{\partial\tilde{\lambda}_i}{\partial\alpha} = \left(\boldsymbol{\Lambda}^{-1}\frac{\partial\textbf{B}}{\partial\alpha}\boldsymbol{\Lambda}\right)_{\left(\xi_i,\xi_i\right)}
\end{equation}

An analytical expression relating the derivative of $\textbf{B}$ with respect to $\alpha$ in function of $\textbf{h}_{z\alpha}$ can be obtained from Equation (\ref{DerivB}). As for fold and BP bifurcations, $\textbf{h}_{z\alpha}$ is computed through finite differences.

From the computational viewpoint, the evaluation of the $\left(2N_H+1\right)n+2$ terms $g_{\alpha}$ represents the main burden of the method when there is a large number of nonlinear DOFs (e.g., in the case of distributed nonlinearities). In this case, parallel computing could help reduce the computational cost
of the algorithm. For localized nonlinearities, one can take advantage of the fact that $\textbf{h}_{z\alpha}$ is a null matrix when $\alpha$ corresponds to Fourier coefficients of linear DOFs.

%%%%%%%%%%%%%%%%%%%%%%%%%%%%
% ---------------- SECTION 3
%%%%%%%%%%%%%%%%%%%%%%%%%%%%

\section{Bifurcation analysis of a satellite structure}

\subsection{Description of the SmallSat spacecraft}

\begin{figure}[ht]
    \centering
    \begin{tabular}{cc}
      \includegraphics[scale=0.8]{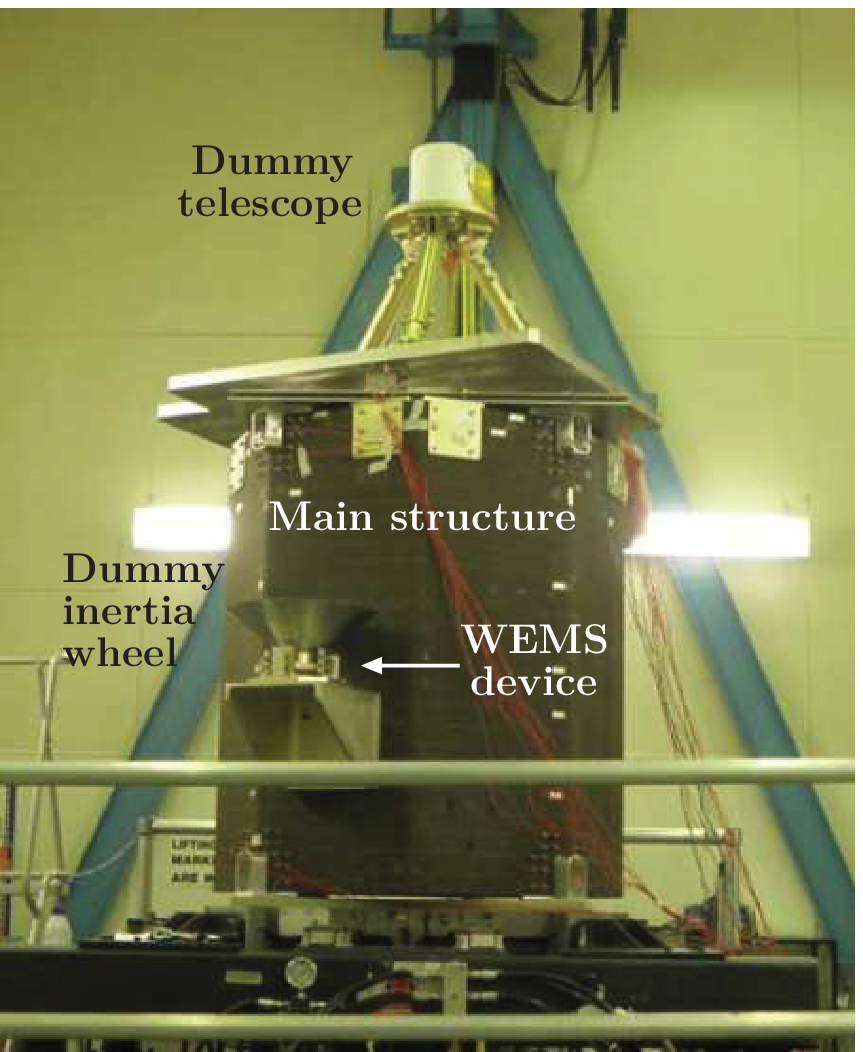} &
      \includegraphics[scale=0.55]{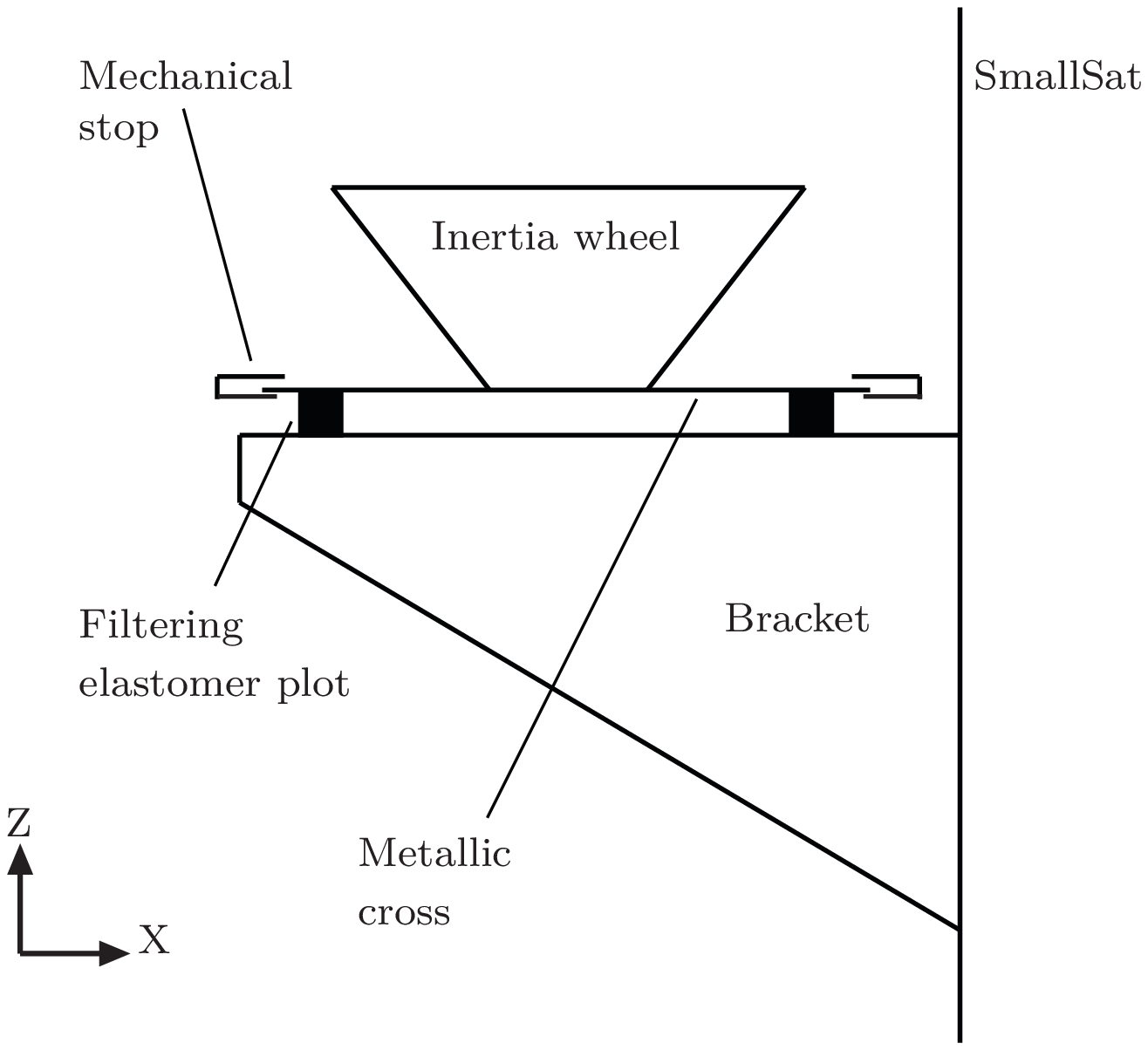} \\ \\
      (a) & (b) \\ \\
    \end{tabular}
    \caption{SmallSat spacecraft. (a) Photograph; (b) Schematic of the WEMS, the nonlinear vibration isolation device. } \label{smallsat_fig}
\end{figure}

The effectiveness of the proposed HB method is demonstrated using the \textit{SmallSat}, a spacecraft structure conceived by EADS-Astrium (now Airbus Defence and Space). The spacecraft is $1.2\,$m in height and $1\,$m in width. A prototype of the spacecraft with a dummy telescope is represented in Figure \ref{smallsat_fig}(a). The nonlinear component, the so-called \textit{wheel elastomer mounting system} (WEMS), is mounted on a bracket connected to the main structure and loaded with an $8$-kg dummy inertia wheel, as depicted in Figure \ref{smallsat_fig}(b). The WEMS device acts a mechanical filter which mitigates the on-orbit micro-vibrations of the inertia wheel through soft elastomer plots located between the metallic cross that supports the inertia wheel and the bracket. To avoid damage of the elastomer plots during launch, the axial and lateral motions of the metallic cross are limited by four nonlinear connections labelled NC1 -- 4. Each NC comprises two mechanical stops that are covered with a thin layer of elastomer to prevent metal-metal impacts.

\begin{figure}[h!t]
    \centering
    \begin{tabular}{c}
      \includegraphics[scale=1.2]{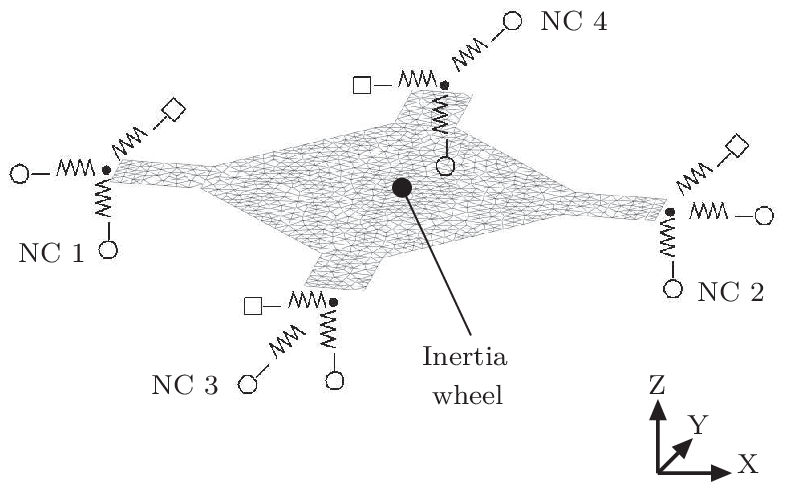} \\ \\
			(a) \\ \\
      \includegraphics[scale=0.65]{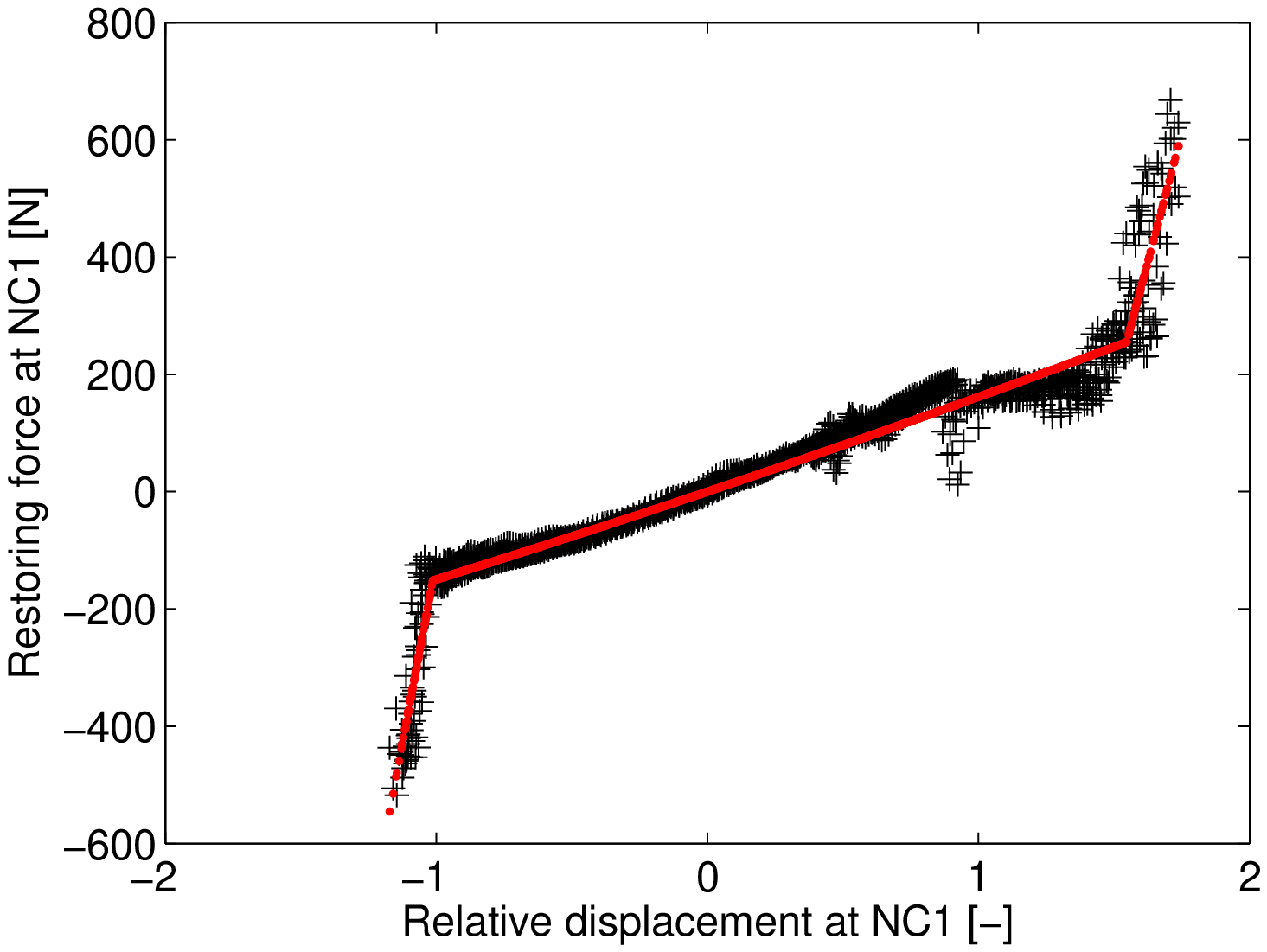} \\ \\
      (b) \\ \\
    \end{tabular}
    \caption{WEMS. (a) Modeling of the device using shell elements, a point mass, linear and nonlinear springs. The linear and nonlinear springs are represented with squares and circles, respectively. (b) Identified stiffness curve of NC1 (in black) and fitting with a trilinear model (in red).} \label{wems_restoring_f}
\end{figure}

A finite element model (FEM) was built to conduct numerical experiments. Linear shell elements were used for the main structure and the instrument baseplate, and a point mass represented the dummy telescope. Proportional damping was considered for these components. As shown in Figure \ref{wems_restoring_f}(a), the metallic cross of the WEMS was also modeled using linear shell elements whereas the inertia wheel was seen as a point mass owing to its important rigidity. To achieve tractable calculations, the linear elements of the FEM were condensed using the Craig-Bampton reduction technique. Specifically, the FEM was reduced to 10 internal modes and 9 nodes (excluding DOFs in rotation), namely both sides of each NC and the inertia wheel. In total, the reduced-order model thus contains 37 DOFs.

Each NC was then modeled using a trilinear spring in the axial direction (elastomer in traction/compression plus two stops), a bilinear spring in the radial direction (elastomer in shear plus one stop) and a linear spring in the third direction (elastomer in shear). The values of the spring coefficients were identified from experimental data \cite{noel2014}. For instance, the stiffness curve identified for NC1 is displayed in Figure \ref{wems_restoring_f}(b). To avoid numerical issues, regularization with third-order polynomials was utilized in the close vicinity of the clearances to implement $C^1$ continuity. The dissipation in the elastomer plots was modeled using lumped dashpots. We note that the predictions of the resulting nonlinear model were found in good agreement with experimental data \cite{noel2014,renson2014}.

For confidentiality, clearances and displacements are given through adimensionalized quantities throughout the paper.

\subsection{Frequency response, Floquet exponents and bifurcation detection}

\begin{figure}[h!t]
    \centering
    \begin{tabular}{c}
      \includegraphics[scale=0.65]{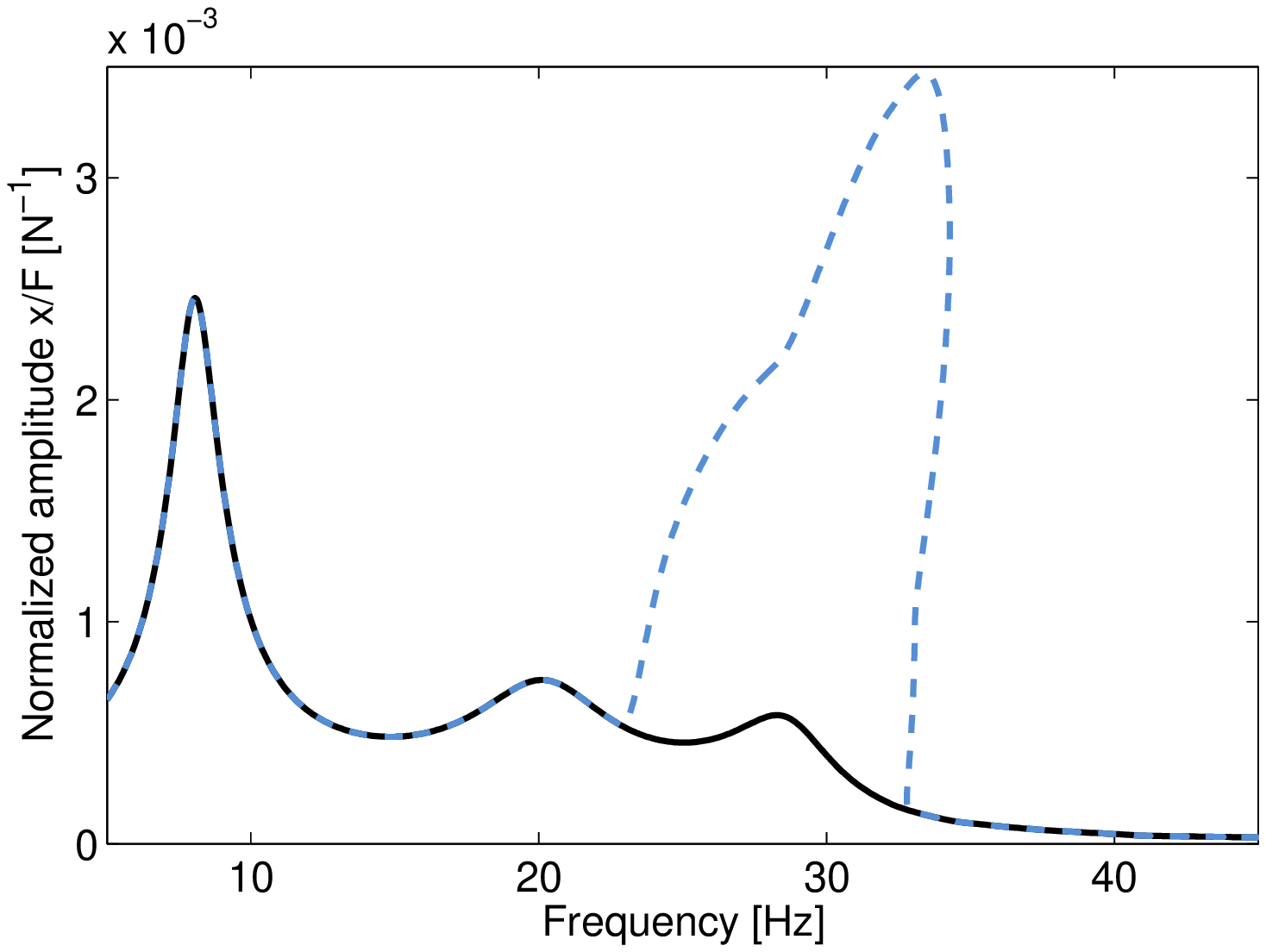} \\ \\
			(a) \\ \\
      \includegraphics[scale=0.65]{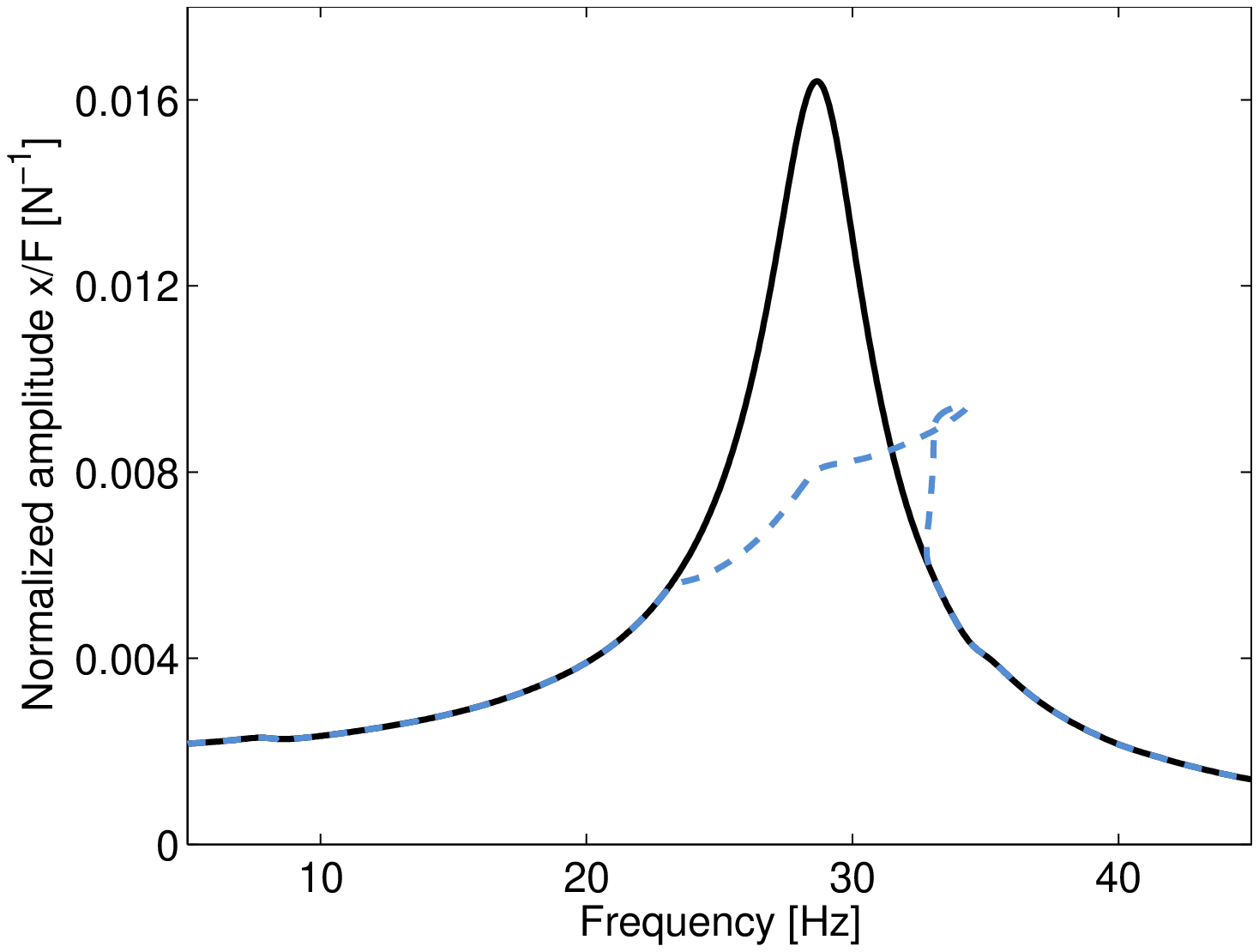} \\ \\
      (b) \\ \\
    \end{tabular}
    \caption{Normalized frequency response at NC1. (a) NC1-$X$; (b) NC1-$Z$. The solid and dashed lines represent forcing amplitude of $F = 50\,$N and $F = 155\,$N, respectively. }\label{low_high_forcing_level}
\end{figure}

The forced response of the satellite for harmonic forcing applied to the vertical DOF of the inertia wheel is computed using HB with $N_H = 9$ harmonics and $N = 1024$ points per period. Figure \ref{low_high_forcing_level} depicts the system's frequency response curves at NC1-$X$ and NC1-$Z$ for two forcing amplitudes, $F = 50$ and $155\,$N. For a clear assessment of the effects of the nonlinearities, the response amplitudes are normalized with the forcing amplitude $F$ in this figure. Because the normalized responses for the two forcing amplitudes coincide up to 23 Hz, one can conclude that the motion is purely linear in this frequency range. Conversely, the mode with a linear resonance frequency of $28.8\,$Hz is greatly affected by the WEMS nonlinearities, which can be explained by the fact that this mode combines bracket deflection with WEMS motion \cite{renson2014}.

\begin{figure}[h!t]
    \centering
    \begin{tabular}{c}
      \includegraphics[scale=0.65]{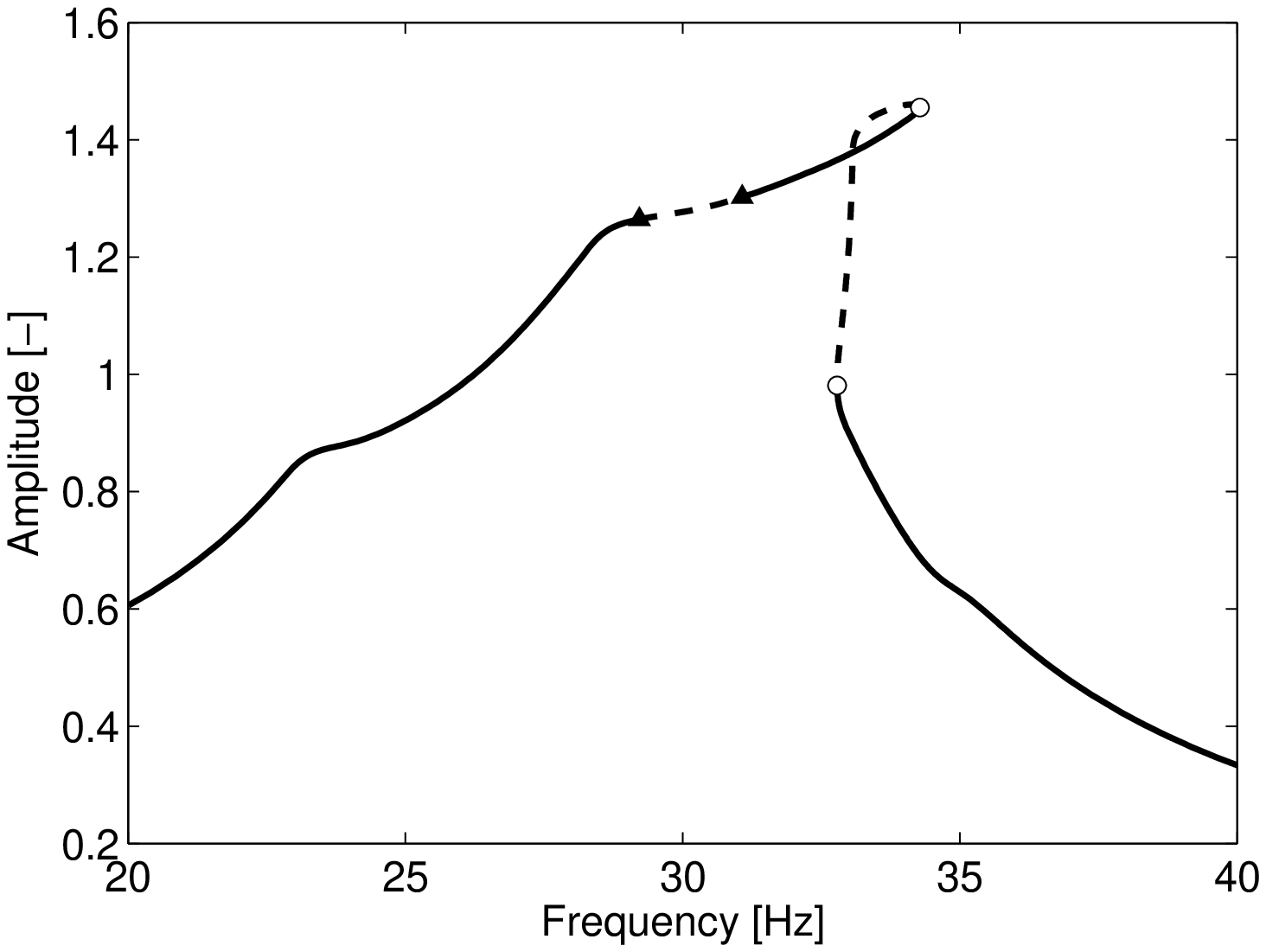} \\ \\
			(a) \\ \\
      \includegraphics[scale=0.65]{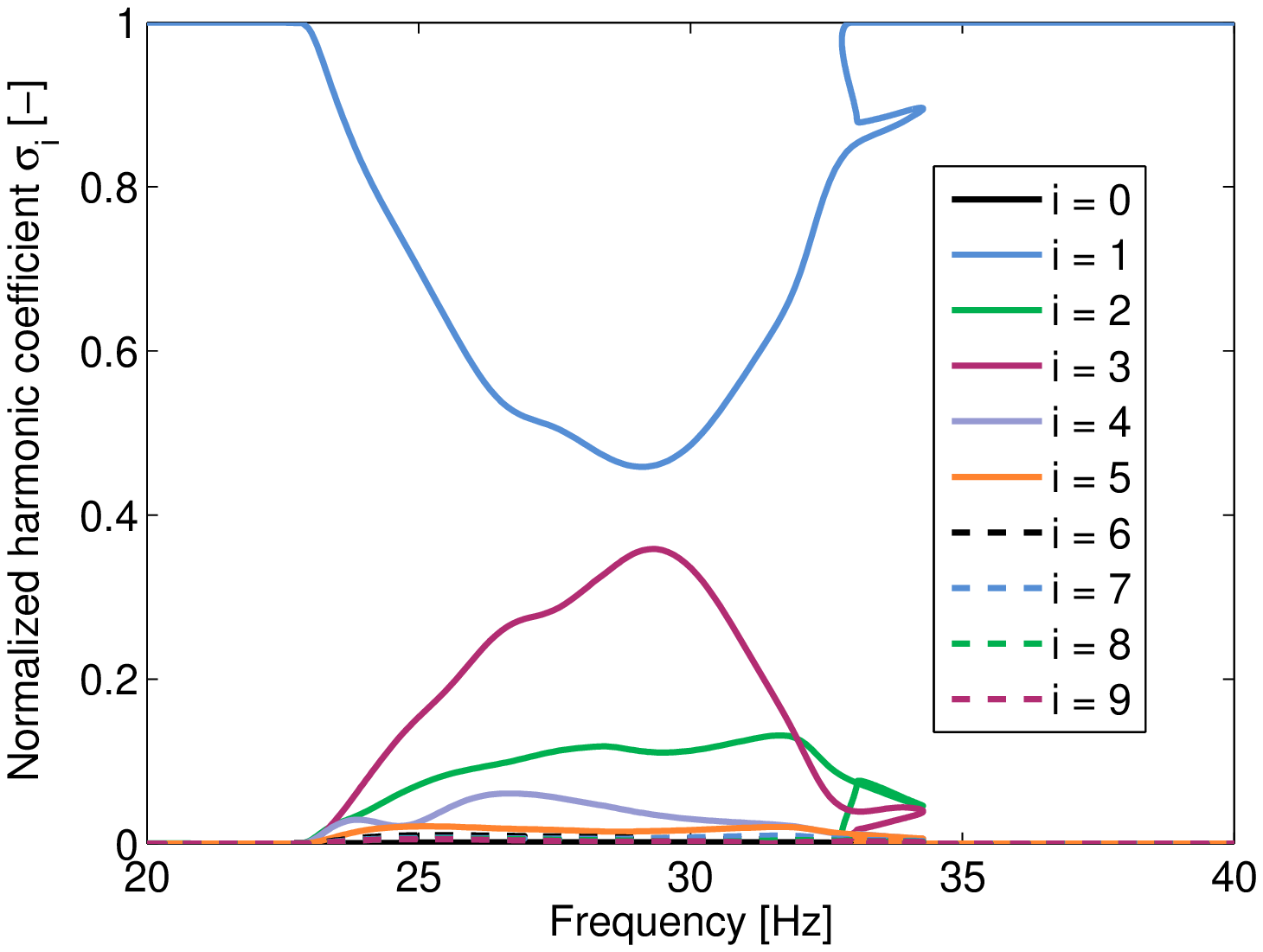} \\ \\
      (b) \\ \\
    \end{tabular}
    \caption{Frequency response at NC1-$Z$ for $F = 155\,$N. (a) Displacement. Circle and triangle markers represent fold and NS bifurcations, respectively. The solid and dashed lines represent stable and unstable solutions, respectively. (b) Harmonic coefficients.} \label{freq_resp_curve}
\end{figure}

Figure \ref{freq_resp_curve}(a) presents a close-up of this nonlinear resonance at NC1-Z where stability and bifurcations are also indicated. The evolution of the normalized harmonic coefficients
\begin{equation}
\sigma_{i} = \frac{\phi_{i}}{\sum_{k = 0}^{N_H}{\phi_{i}}} \ \ \left(i = 0,\hdots,N_H\right)
\end{equation}
with
\begin{equation}
\phi_{0} = \frac{c_0^x}{\sqrt{2}},\quad \phi_{i} = \sqrt{\left(s_i^x\right)^2+\left(c_i^x\right)^2} \ \ \left(i = 1,\hdots,N_H\right)
\end{equation}
is shown in Figure \ref{freq_resp_curve}(b). From $20$ to $23\,$Hz, only the fundamental harmonic is present in the response. In the resonance region, the SmallSat nonlinearities activate other harmonics in the response. Even harmonics contribute to the dynamics because of the asymmetric modeling of the NCs of the WEMS. From the figure, it is also clear that the 6th and higher harmonics have a negligible participation in the response; for this reason, $N_H = 5$ is considered throughout the rest of the paper. 

\begin{figure}[h!t]
    \centering
    \begin{tabular}{cc}
      \includegraphics[scale=0.65]{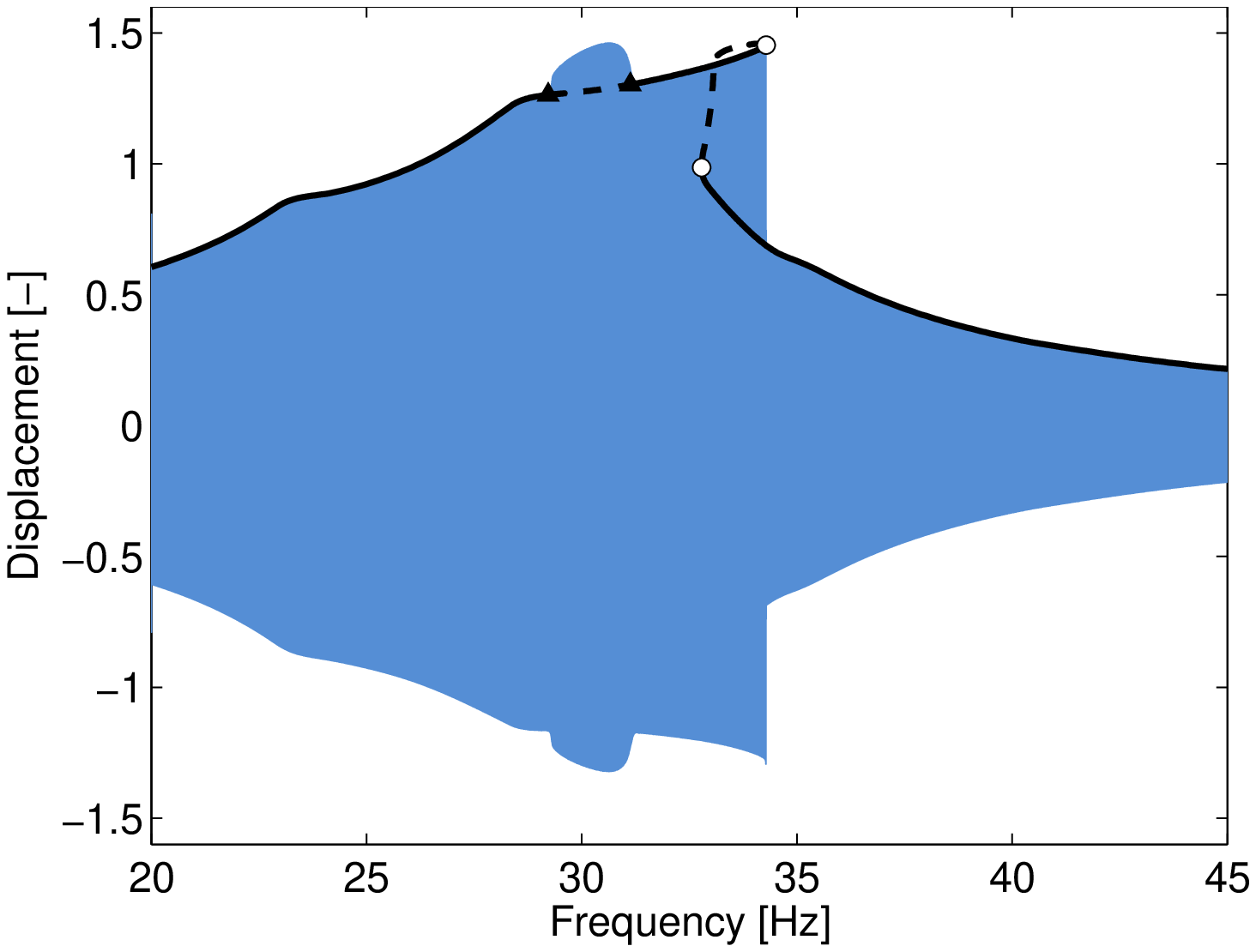} & \includegraphics[scale=0.65]{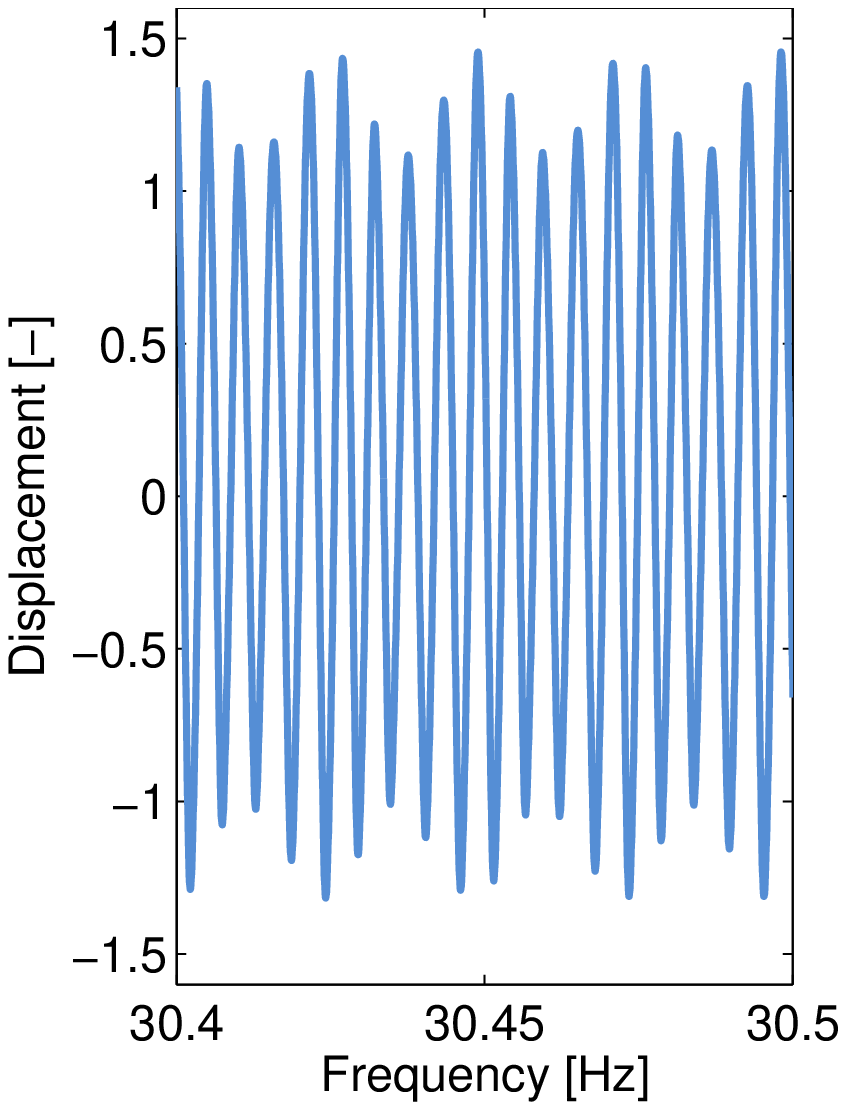} \\ \\
			(a) & (b)\\ \\
      \multicolumn{2}{c}{\includegraphics[scale=0.65]{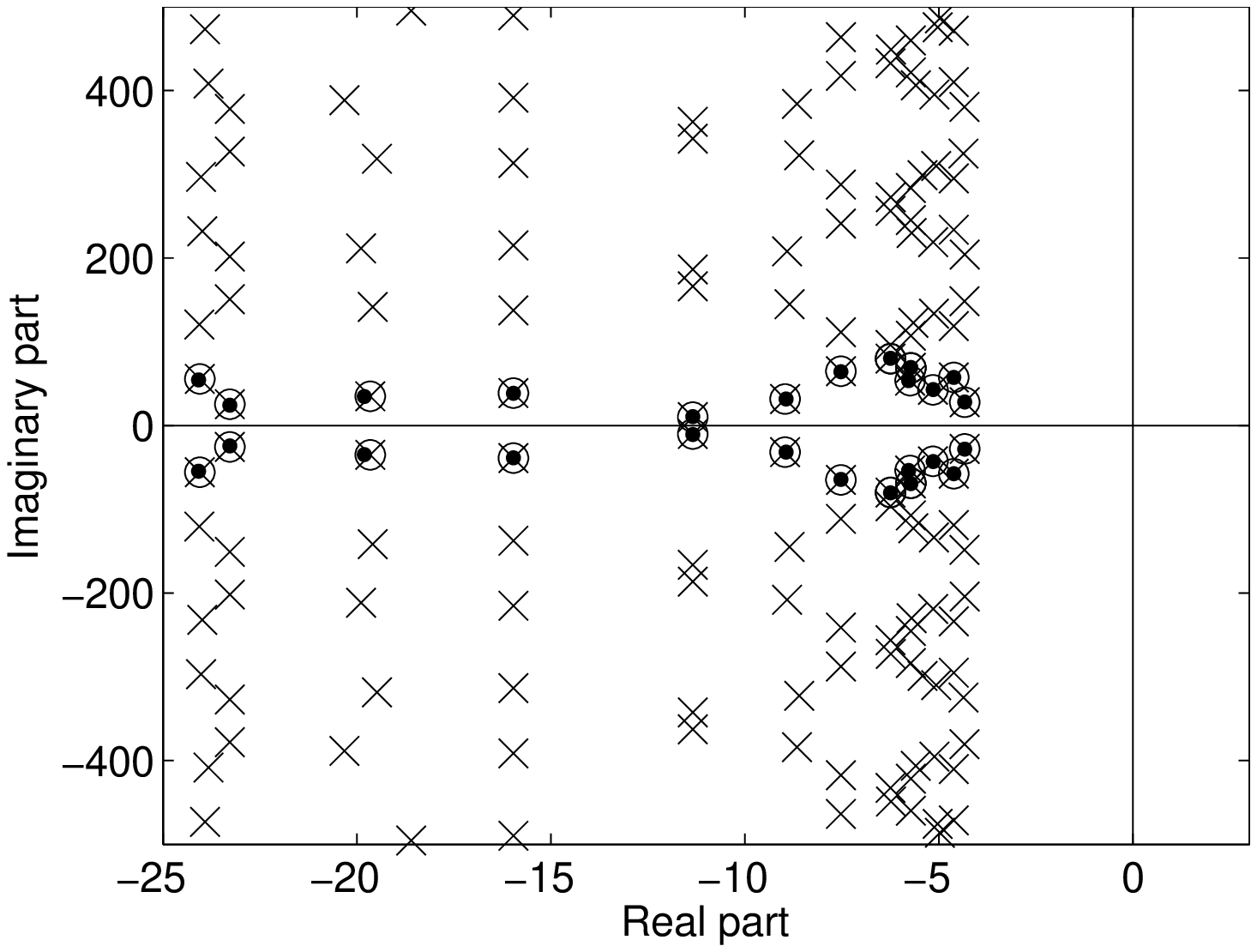}} \\ \\
      \multicolumn{2}{c}{(c)} \\ \\
    \end{tabular}
    \caption{Bifurcations, quasiperiodic oscillations and Floquet exponents. (a) Comparison between HB (black) and a swept-sine response calculated using Newmark's algorithm (blue) for $F = 155\,$N at NC1-$Z$. Circle and triangle markers depict fold and NS bifurcations, respectively. The solid and dashed lines represent stable and unstable solutions, respectively. (b) Close-up of the quasiperiodic oscillations. (c) Periodic solution at $28\,$Hz: Hill's coefficients $\boldsymbol{\lambda}$ (crosses), Floquet exponents $\tilde{\boldsymbol{\lambda}}$ obtained with Hill's method (circles), and Floquet exponents $\tilde{\boldsymbol{\lambda}}_{TI}$ obtained from the monodromy matrix (dots).}\label{freq_resp_comparison}
\end{figure}

The nonlinear resonance has a complex and rich topology. Not only the frequency response curve was found to fold backwards, but bifurcations were detected along the branch by the HB algorithm in Figure \ref{freq_resp_curve}(a). Fold bifurcations are generic for nonlinear resonances and are responsible for stability changes. NS bifurcations imply that further investigation of the dynamics should be carried out in the corresponding range of frequencies, because a branch of quasiperiodic solutions bifurcates out from the main branch. To this end, the response to a swept-sine excitation with a forcing amplitude $F = 155\,$N and with a sweep rate of $0.5\,$Hz/min was computed with a Newmark time integration scheme. The sampling frequency was chosen very high, i.e., $3000\,$Hz, to guarantee the accuracy of the simulation in Figure \ref{freq_resp_comparison}(a). The first observation is that overall the frequency response computed using HB provides a very accurate estimation of the envelope of the swept-sine response. In addition, a modulation of the displacement's envelope is clearly noticed between the two NS bifurcations highlighted by HB, and the close-up in Figure \ref{freq_resp_comparison}(b) confirms the presence of quasiperiodic oscillations. Interestingly, the amplitudes associated with these oscillations are slightly larger than those at resonance, which shows the importance of a proper characterization of these oscillations.

Figure \ref{freq_resp_comparison}(c) shows a subset of Hill's coefficients $\boldsymbol{\lambda}$ and Floquet exponents $\tilde{\boldsymbol{\lambda}}$ obtained with Hill's method ($N_H = 5$ and $N = 1024$) for the stable periodic solution at $28\,$Hz in Figure \ref{freq_resp_comparison}(a). Floquet exponents $\tilde{\boldsymbol{\lambda}}_{TI}$ are also calculated from the monodromy matrix evaluated with a Newmark time integration scheme as in \cite{peeters2009}; they serve as a reference solution. The comparison demonstrates that the actual Floquet exponents corresponds to the Hill's coefficients that are the closest ones to the real axis, which validates the sorting criterion discussed in Section \ref{stab_anal_sec}. The other coefficients are spurious, but they seem to be aligned according to the same pattern as that of the actual Floquet exponents. The location of all exponents in the left-half plane indicates that the solution at $28\,$Hz is indeed stable.

\begin{figure}[h!t]
    \centering
    \begin{tabular}{cc}
      \includegraphics[scale=0.65]{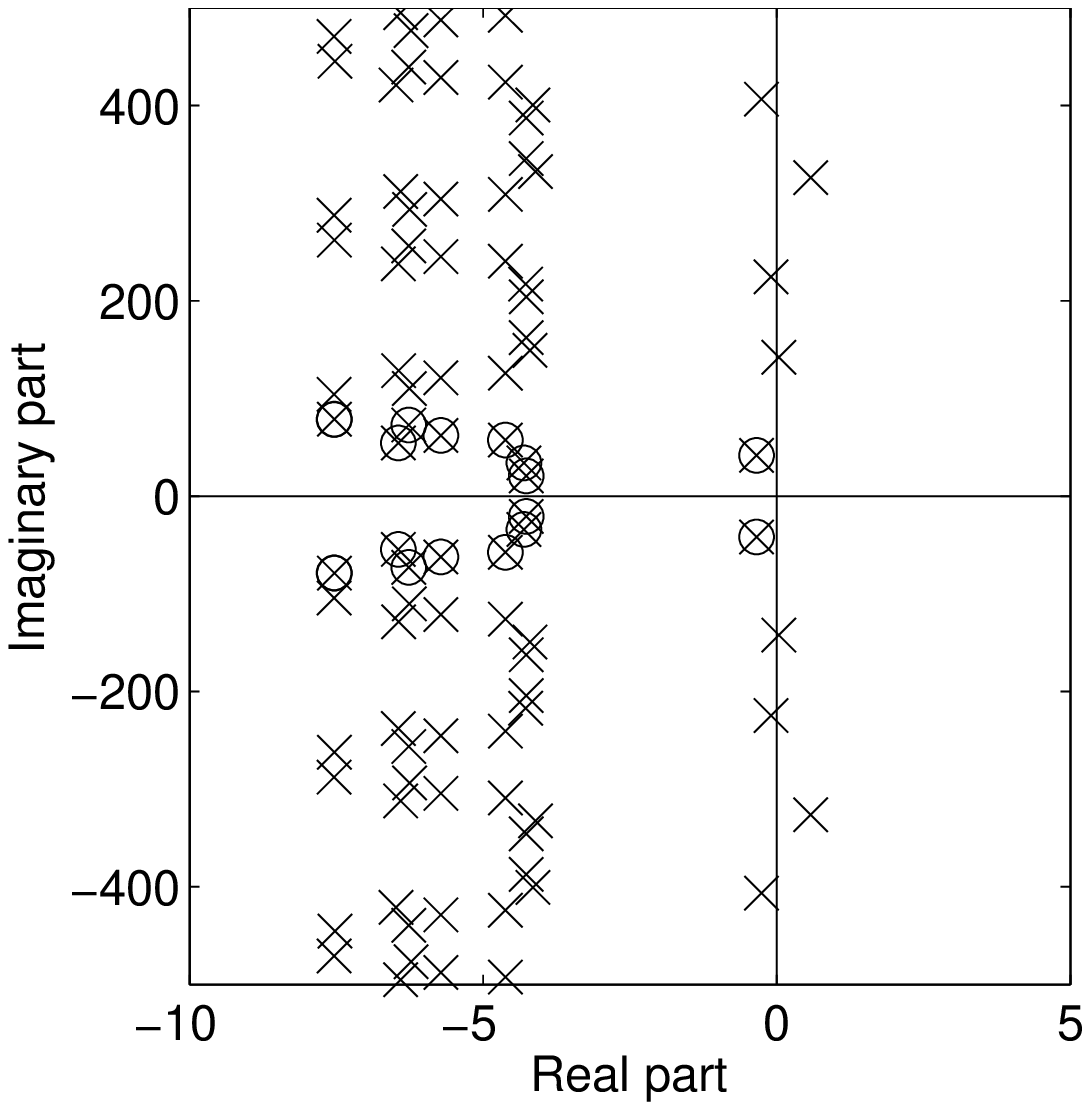} & \includegraphics[scale=0.65]{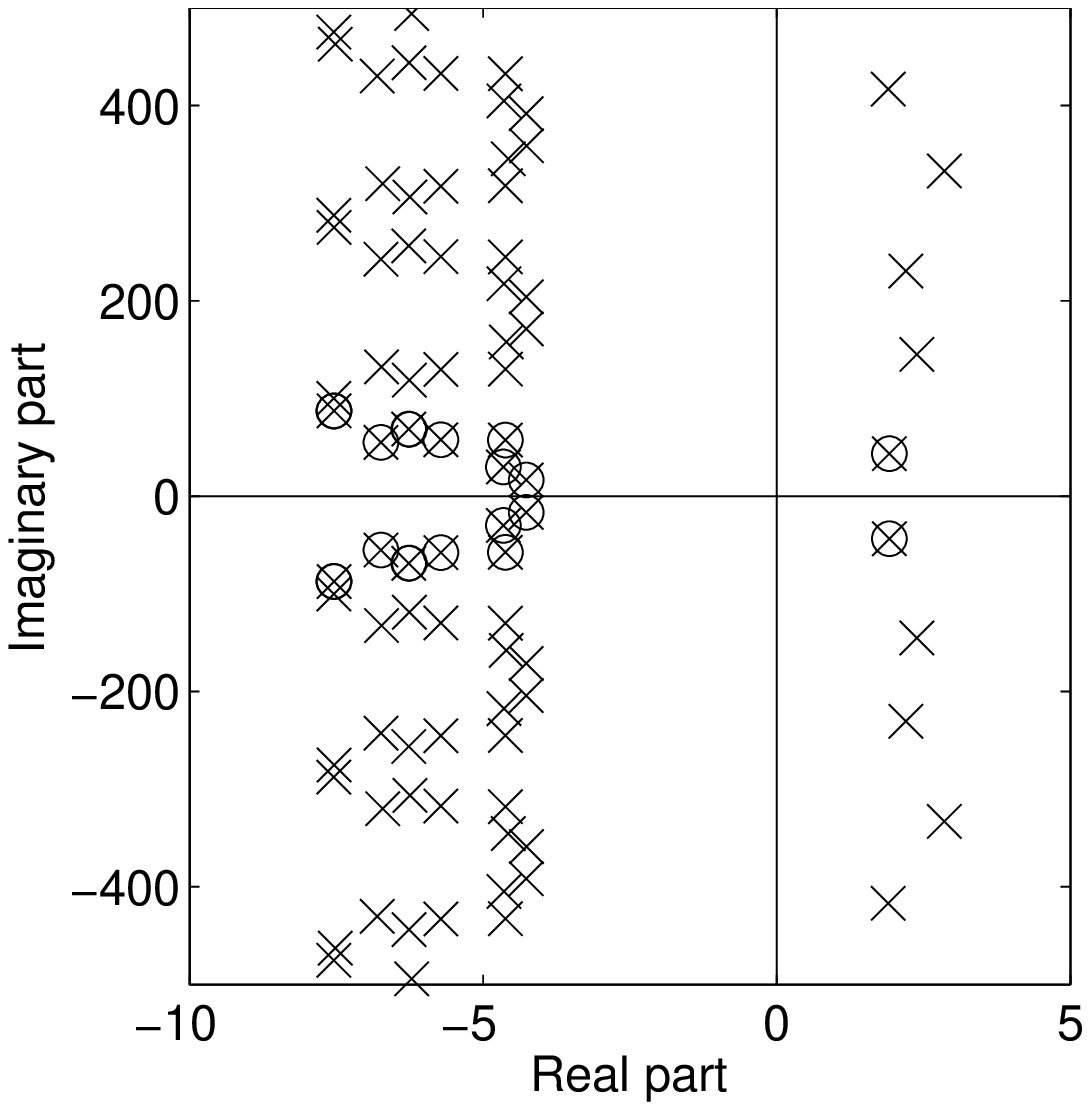} \\ \\
			(a) & (b)\\ \\
      \includegraphics[scale=0.65]{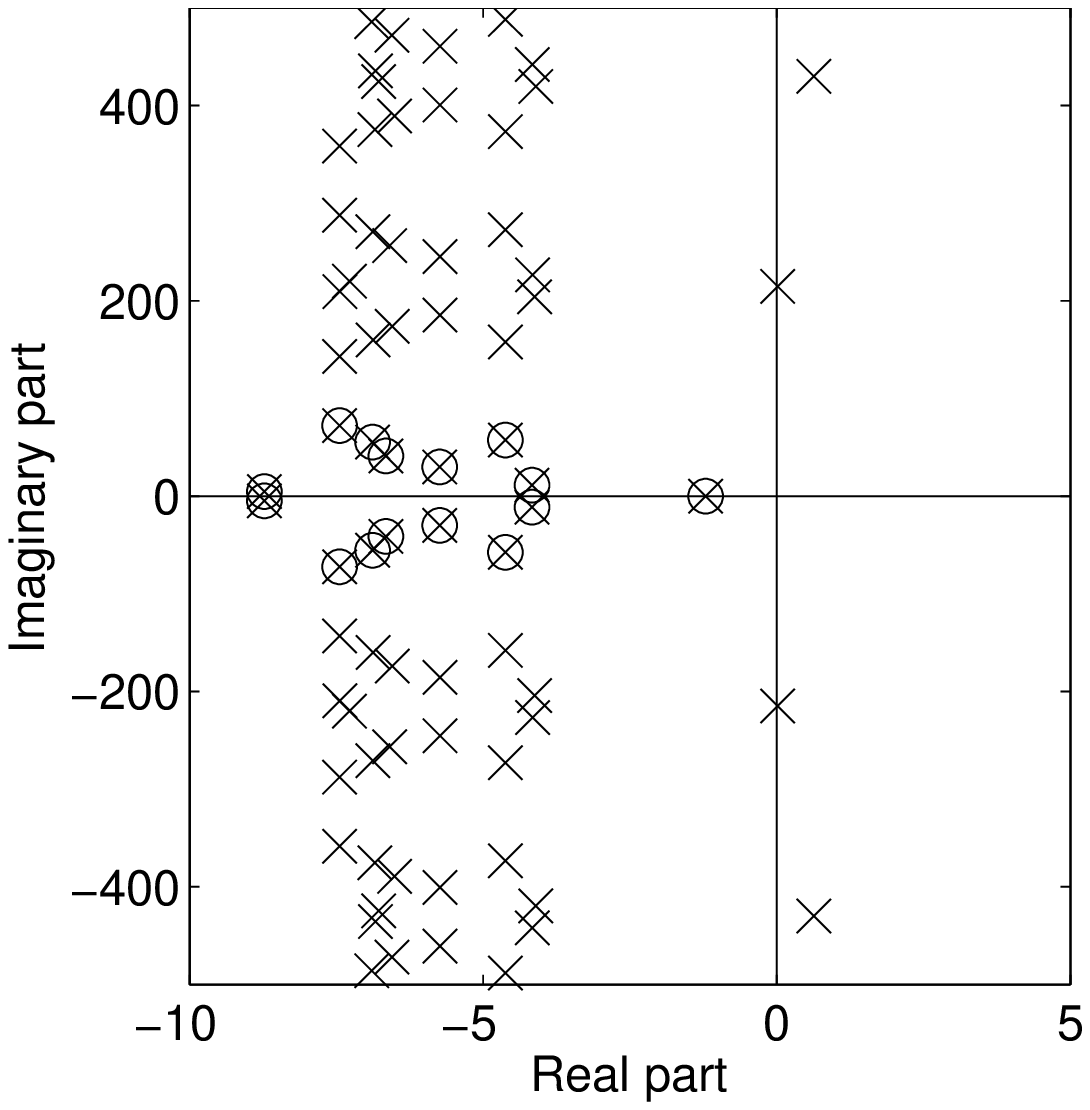} & \includegraphics[scale=0.65]{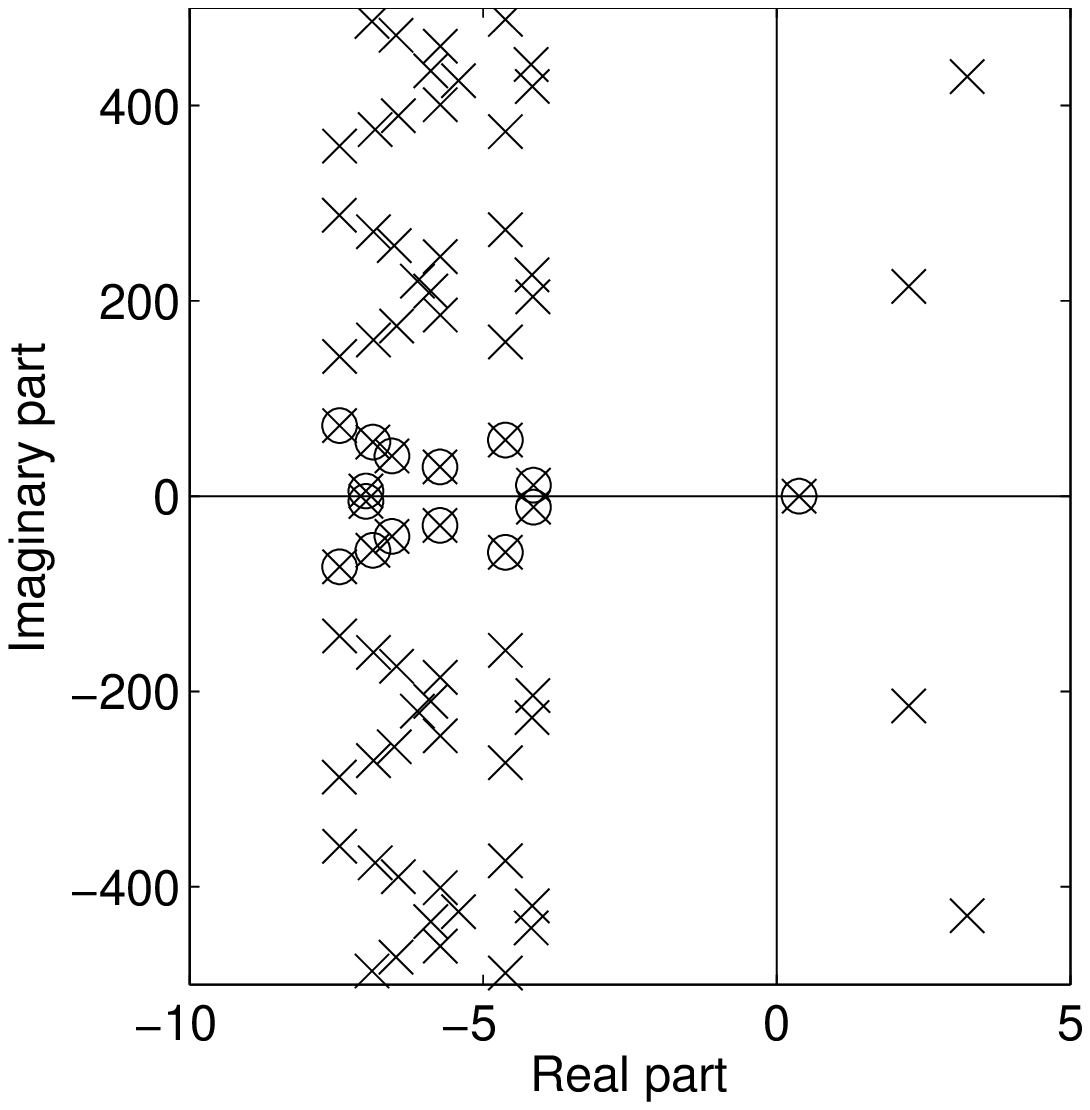}\\ \\
      (c) & (d) \\ \\
    \end{tabular}
    \caption{Floquet exponents in the vicinity of the fold and NS bifurcations. (a) Before the NS bifurcation at $\omega = 29.11\,$Hz, stable region. (b) After the NS bifurcation at $\omega = 29.8\,$Hz, unstable region. (c) Before the fold bifurcation at $\omega = 34.25\,$Hz, stable region. (d) After the fold bifurcation at $\omega = 34.28\,$Hz, unstable region. Hill's coefficients $\boldsymbol{\lambda}$ and Floquet exponents $\tilde{\boldsymbol{\lambda}}$ obtained with Hill's method are denoted with crosses and circles, respectively.}\label{floquet_comp}
\end{figure}

The evolution of Hill's coefficients and Floquet exponents in the vicinity of the first NS bifurcation is given in Figures \ref{floquet_comp}(a-b). Before the bifurcation, the Floquet exponents lie all in the left-half plane in Figure \ref{floquet_comp}(a), which indicates a stable solution. We stress that a pair of Hill's coefficients has already crossed the imaginary axis in this figure, which evidences that considering all Hill's coefficients can thus lead to misjudgement. After the bifurcation, a pair of complex conjugates Floquet exponents now lie in the right-half plane in Figure \ref{floquet_comp}(b), which means that the system underwent a NS bifurcation and lost stability. A similar scenario is depicted for the first fold bifurcation in Figures \ref{floquet_comp}(c-d) with the difference that a single Floquet exponent crosses the imaginary axis through zero.

\subsection{Bifurcation tracking}

\begin{figure}[h!t]
    \centering
		\begin{tabular}{c}
      \includegraphics[scale=0.65]{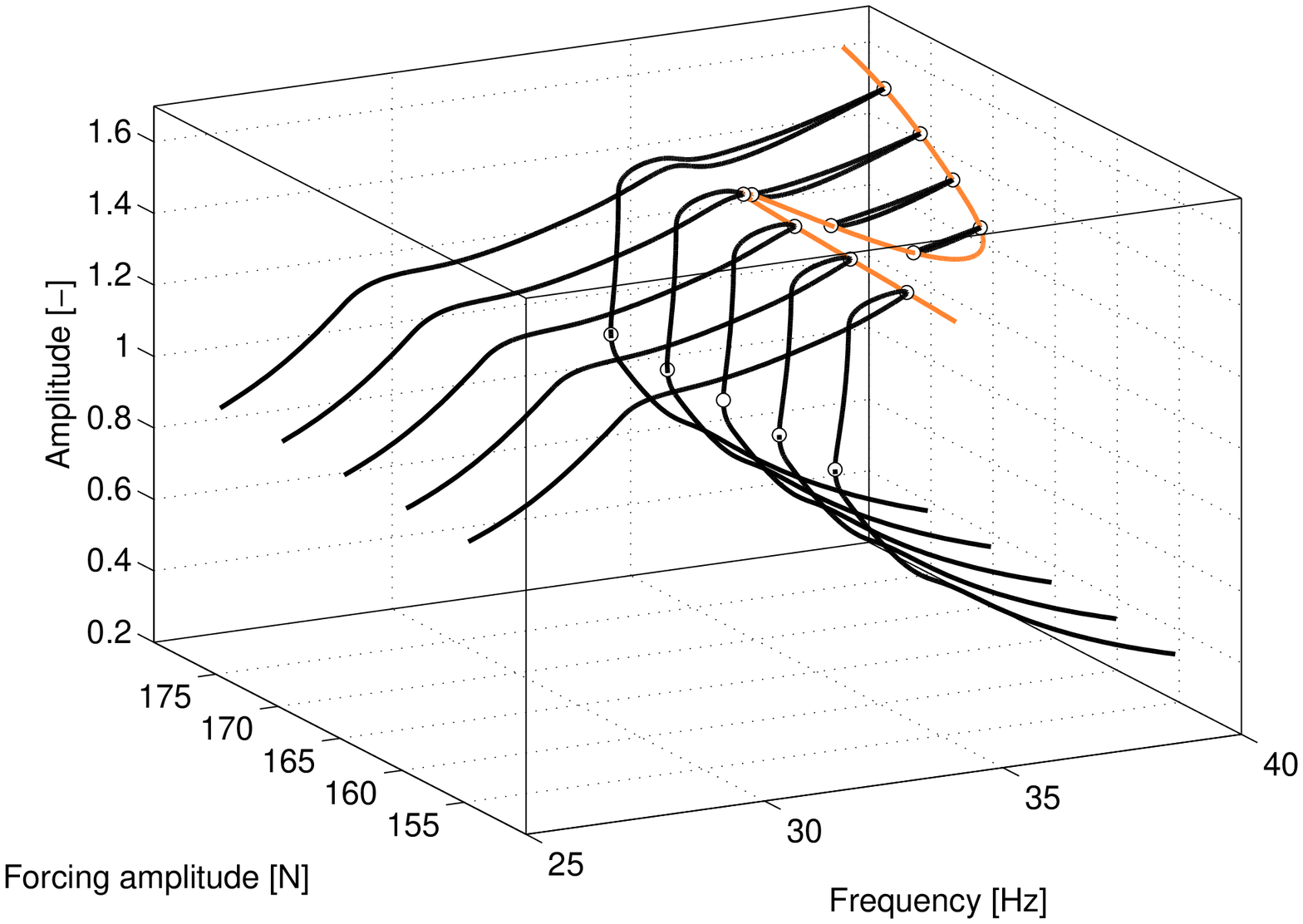} \\ \\
			(a) \\ \\
			\includegraphics[scale=0.65]{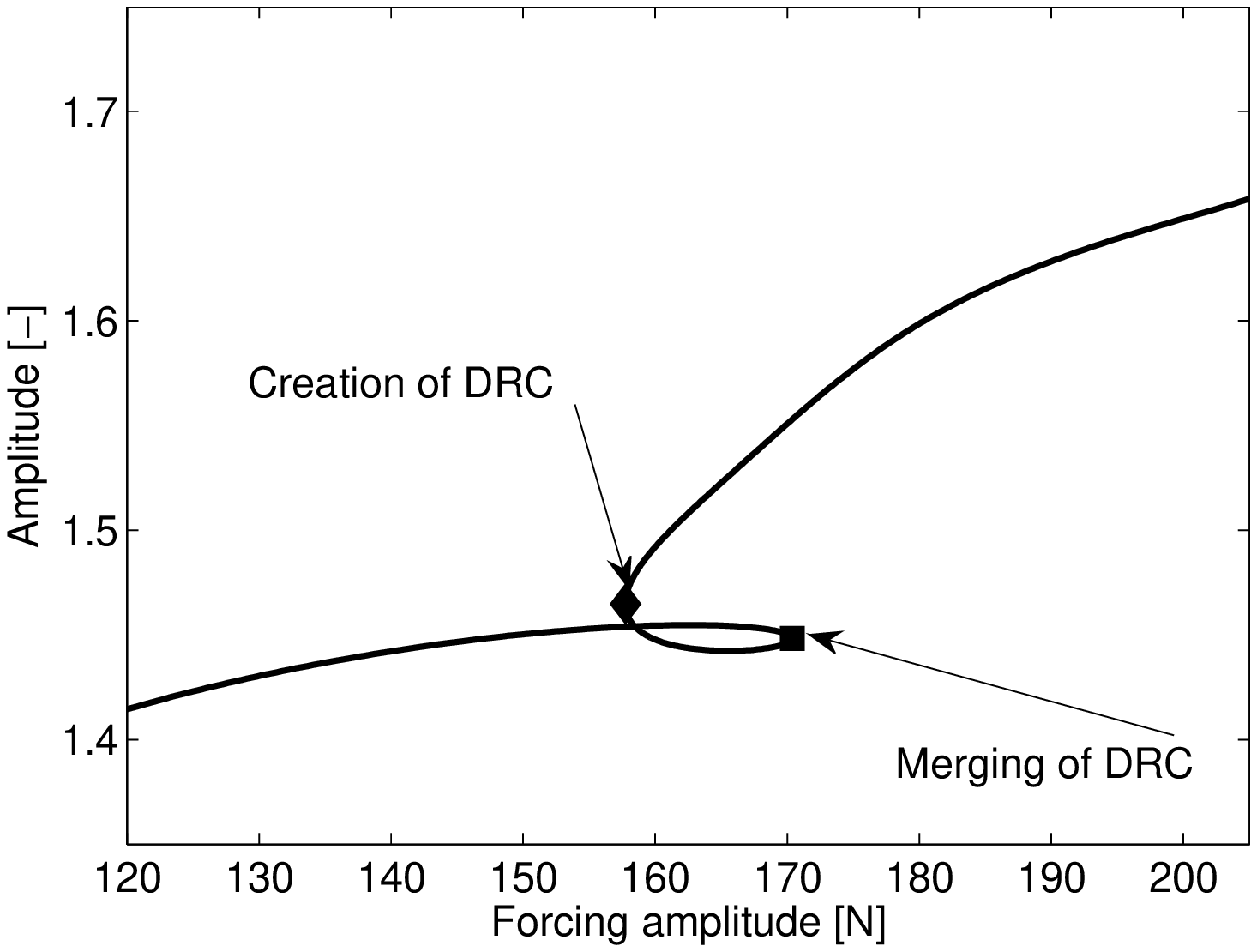} \\ \\
      (b) \\ \\
			\end{tabular}
    \caption{Tracking of the fold bifurcations of the resonance peak. (a) Three-dimensional space. Orange line: branch of fold bifurcations; black lines: frequency responses at NC1-$Z$ for $F = 155\,$N, $160\,$N, $170\,$N and $175\,$N. Circle markers depict fold bifurcations. (b) Two-dimensional projection of the branch of fold bifurcations.} \label{LP_tracking}
\end{figure}

\begin{figure}[h!t]
\centering \includegraphics[scale=0.65]{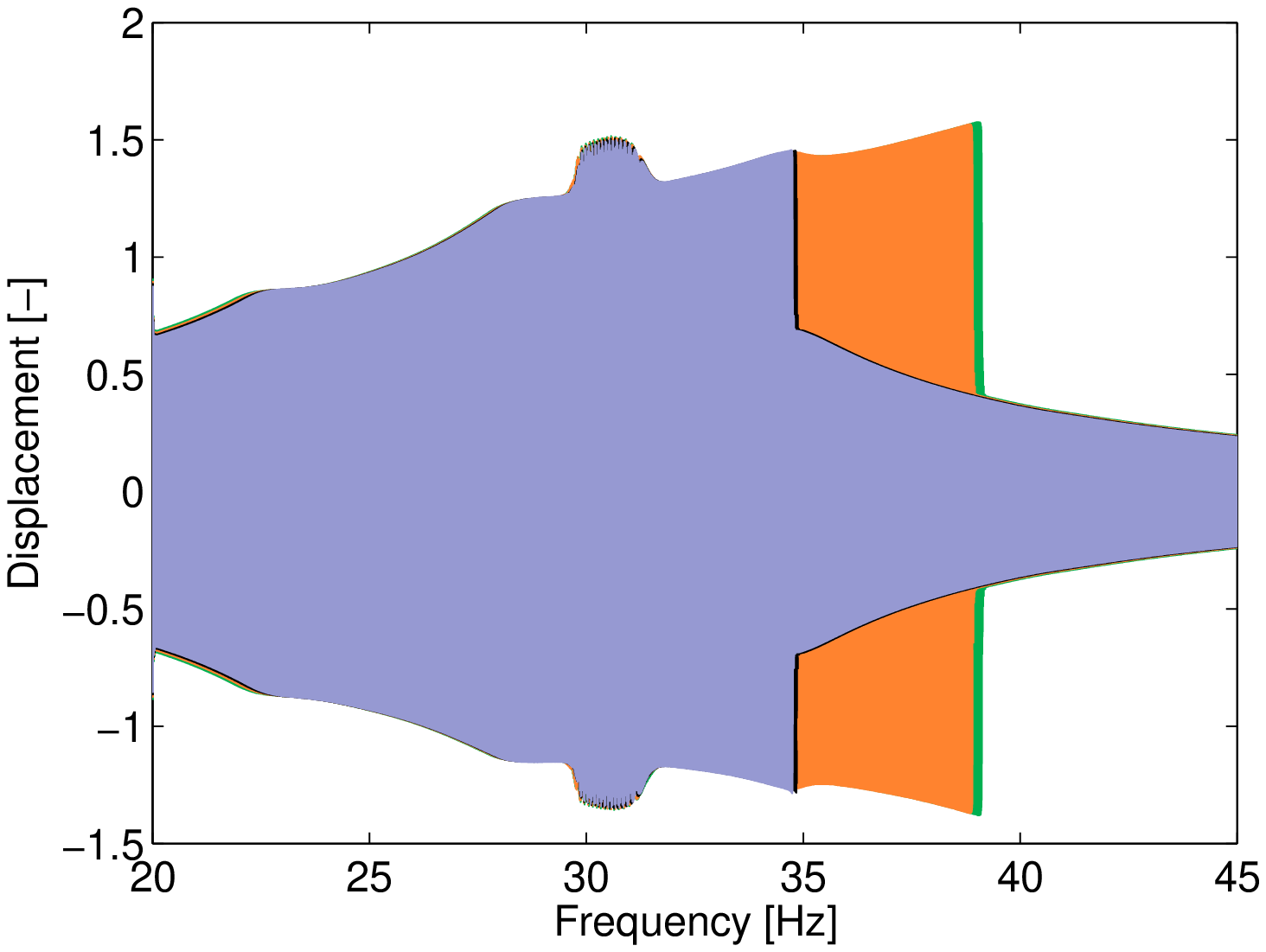}
\caption{NC1-$Z$ displacement for swept-sine excitations with $F = 168\,$N (purple), $170\,$N (black), $172\,$N (orange) and $174\,$N (green).} \label{jump_phenomenon}
\end{figure}

The fold bifurcations revealed in the previous section are now tracked in the codimension-2 forcing frequency-forcing amplitude space using the algorithm presented in Section \ref{biftracking}. Figure \ref{LP_tracking}(a) represents the resulting fold curve, together with the frequency responses of the system for different forcing amplitudes. Very interestingly, the algorithm initially tracks the fold bifurcations of the main frequency response, but it then turns back to reveal a \textit{detached resonance curve} (DRC), or \textit{isola}. Such an attractor is rarely observed for real structures in the literature. Figure \ref{LP_tracking}(b), which shows the projection of the fold curve in the forcing amplitude-response amplitude plane, highlights that the DRC is created when $F = 158\,$N. The DRC then expands both in frequency and amplitude until $F = 170\,$N for which merging with the resonance peak occurs. Because the upper part of the DRC is stable, the resonance peak after the merging is characterized by a greater frequency and amplitude. This merging process is further illustrated in Figure \ref{jump_phenomenon}, where the responses to swept-sine excitations of amplitude $168$, $170$, $172$ and $174\,$N are superimposed. Moving from $F = 170\,$N to $F = 172\,$N leads to a sudden `jump' in resonance frequency and amplitude.

\begin{figure}[h!t]
    \centering
		\begin{tabular}{c}
      \includegraphics[scale=0.65]{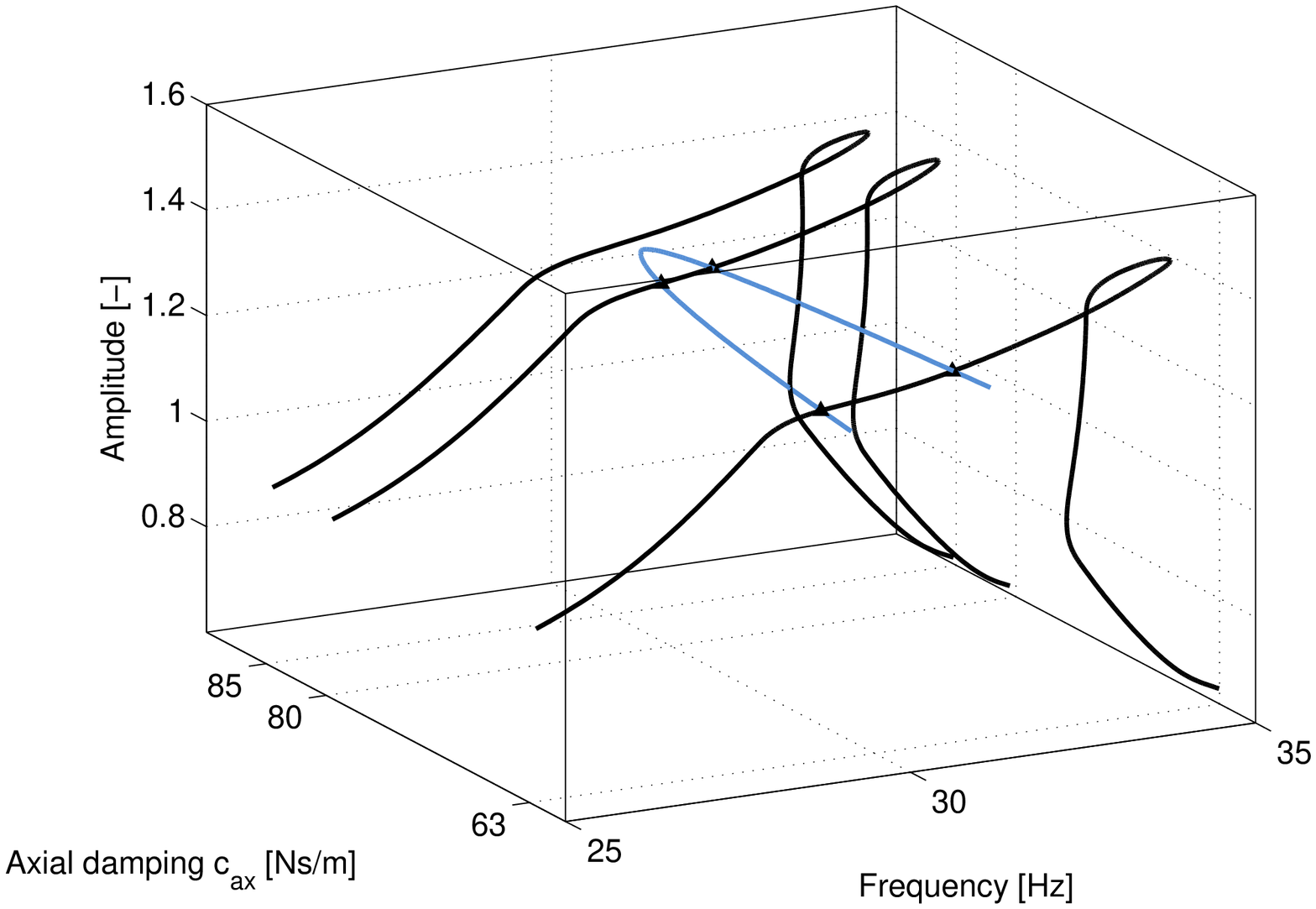} \\ \\
			(a) \\ \\
			\includegraphics[scale=0.65]{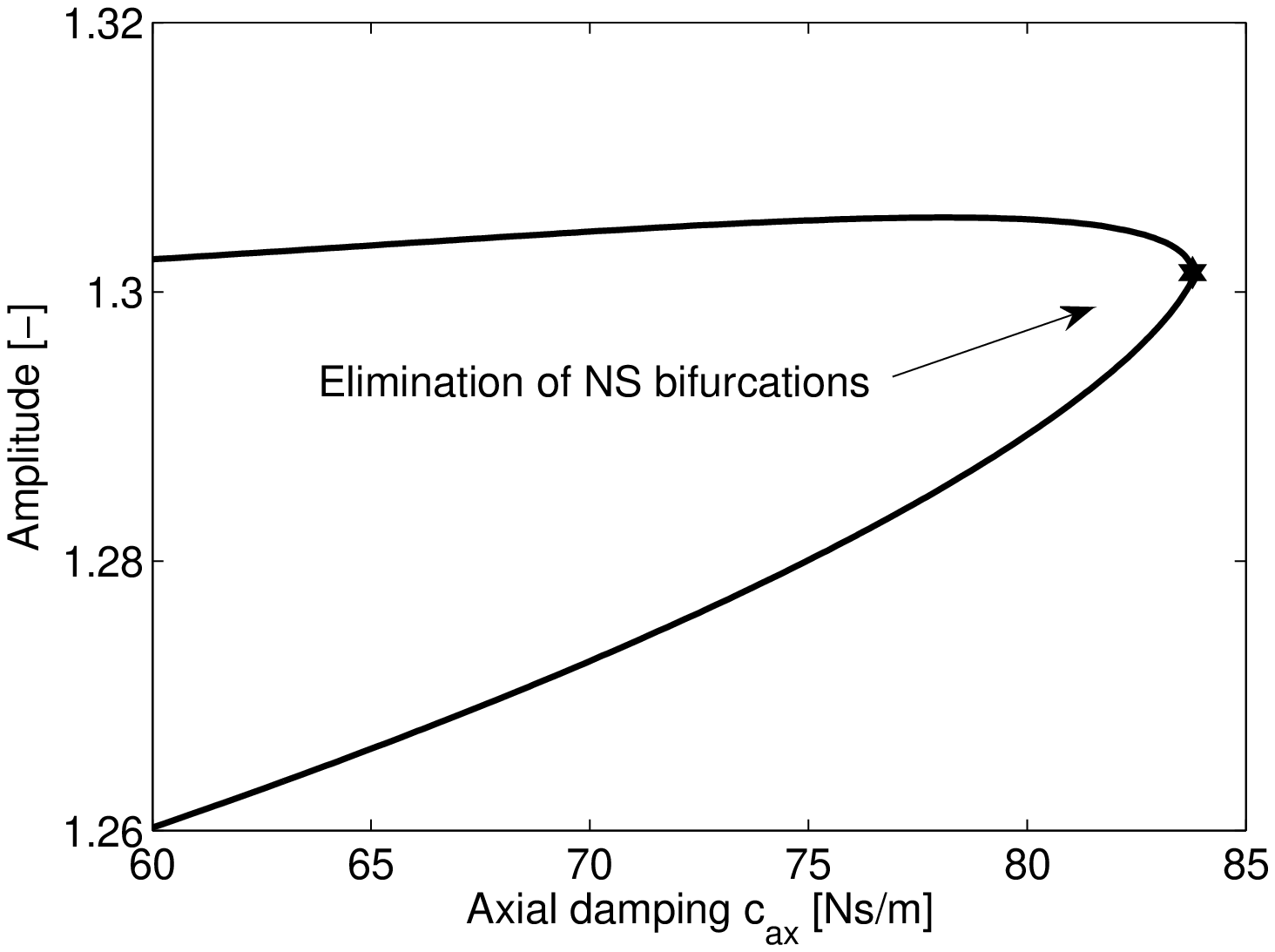} \\ \\
      (b) \\ \\
			\end{tabular}
    \caption{Tracking of the NS bifurcations. (a) Three-dimensional space. Blue line: branch of NS bifurcations; black lines: frequency responses of the at NC1-$Z$ for $F = 155\,$N and for $c_{ax} = 63\,$Ns/m, $80\,$Ns/m and $85\,$Ns/m. Triangle markers depict NS bifurcations. (b) Two-dimensional projection of the branch of NS bifurcations.} \label{NS_tracking}
\end{figure}

Bifurcation tracking is not only useful for understanding the system's dynamics and revealing additional attractors, but it can also be used for engineering design. For instance, we use it herein to study the influence of the axial dashpot $c_{ax}$ of the WEMS device on the quasiperiodic oscillations. Figure \ref{NS_tracking}(a) depicts the NS curve in the codimension-2 forcing frequency-axial damping space to which frequency responses computed for $c_{ax} = 63\,$Ns/m (reference), $80\,$Ns/m and $85\,$Ns/m are superimposed. The projection in the $c_{ax}$-response amplitude plane in Figure \ref{NS_tracking}(b) shows that the two NS bifurcations, and hence quasiperiodic oscillations, can be completely eliminated for $c_{ax}=84\,$Ns/m.

\begin{figure}[h!t]
    \centering
		\begin{tabular}{c}
      \includegraphics[scale=0.65]{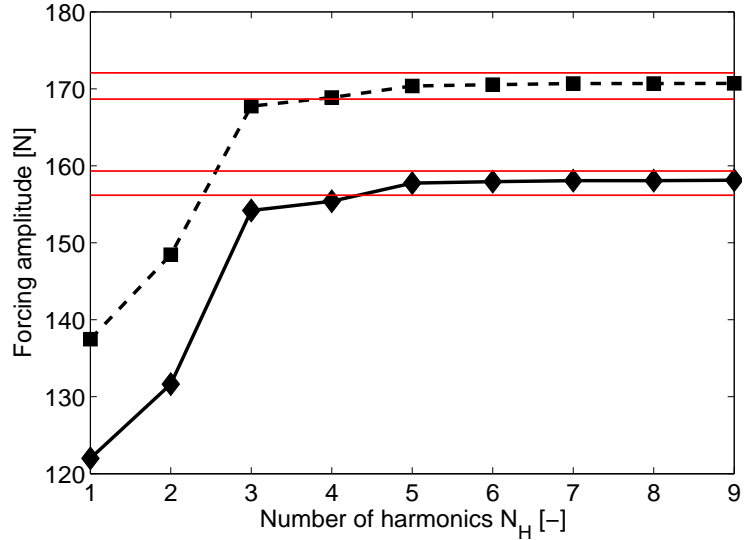} \\ \\
			(a) \\ \\
			\includegraphics[scale=0.65]{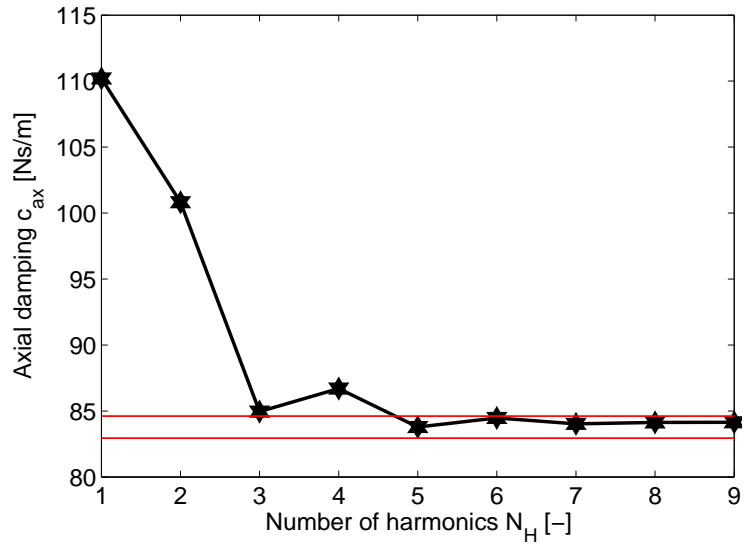} \\ \\
      (b) \\ \\
			\end{tabular}
    \caption{Convergence of the bifurcation tracking algorithm with respect to the number of harmonics $N_H$. (a) Fold bifurcations: the solid line with diamond markers corresponds to the creation of the DRC, whereas the dashed line with square markers corresponds to the merging of the DRC with the main resonance. (b) NS bifurcations: the solid line with stars markers corresponds to the elimination of the bifurcations. For both figures, the red lines gives $+1\%$ and $-1\%$ variations with respect to $N_H$ = 5.} \label{LP_NS_convergence}
\end{figure}

Finally, the convergence of the algorithm is assessed in Figure \ref{LP_NS_convergence}, for which the fold and NS curves presented in Figures \ref{LP_tracking} and \ref{NS_tracking} were recomputed for a number of harmonics $N_H$ varying from $1$ to $9$. In both cases, a clear convergence of the results is observed. These figures also confirm that retaining $5$ harmonics in the Fourier series leads to an error less than 1\%.
 
\FloatBarrier
\section{Conclusions}

The purpose of this paper was to exploit the harmonic balance method for the computation of periodic solutions and their bifurcations in codimension-2 parameter space. Particular attention was devoted to the computational efficiency of the algorithm, which motivated the use of a bordering technique and the development of a new procedure for Neimark-Sacker bifurcations exploiting the properties of eigenvalue derivatives.

The method was demonstrated using a real satellite example possessing several mechanical stops. By combining the proposed algorithm with Craig-Bampton reduction, we showed that bifurcation tracking of large-scale structures with localized nonlinearities is now within reach; it represents a particularly effective tool both for uncovering dynamical attractors and for engineering design.

\section*{Acknowledgments}

\medskip
The authors Thibaut Detroux, Luc Masset and Gaetan Kerschen would like to acknowledge the financial support of the European Union (ERC Starting Grant NoVib 307265). The author L. Renson is a Marie-Curie COFUND Postdoctoral Fellow of the University of Liège, co-funded by the European Union.

\bibliography{HB_paper_V8}

\begin{thebibliography}{10}

\bibitem{padmanabhan1995}
C.~Padmanabhan and R.~Singh, ``Analysis of periodically excited non-linear
  systems by a parametric continuation technique,'' {\em Journal of Sound and
  Vibration}, vol.~184, no.~1, pp.~35--58, 1995.

\bibitem{stoykov}
S.~Stoykov and S.~Margenov, ``Numerical computation of periodic responses of
  nonlinear large-scale systems by shooting method,'' {\em Computers and
  Mathematics with Applications}, vol.~367, no.~5, pp.~2257--2267, 2014.

\bibitem{peeters2009}
M.~Peeters, R.~Vigui{\'e}, G.~S{\'e}randour, G.~Kerschen, and J.-C. Golinval,
  ``Nonlinear normal modes, {Part II}: Toward a practical computation using
  numerical continuation techniques,'' {\em Mechanical Systems and Signal
  Processing}, vol.~23, no.~1, pp.~195--216, 2009.

\bibitem{doedel1997}
E.~J. Doedel, A.~R. Champneys, T.~F. Fairgrieve, Y.~A. Kuznetsov, B.~Sandstede,
  and X.~Wang, ``\textsc{auto97}: Continuation and bifurcation software for
  ordinary differential equations (with homcont),'' {\em User’s Guide,
  Concordia University, Montreal, Canada, Available from
  http://indy.cs.concordia.ca}, 1997.

\bibitem{ascher1979}
U.~Ascher, J.~Christiansen, and R.~D. Russell, ``A collocation solver for mixed
  order systems of boundary value problems,'' {\em Mathematics of Computation},
  pp.~659--679, 1979.

\bibitem{kuznetsov1995}
Y.~A. Kuznetsov and V.~V. Levitin, ``\textsc{content}: A multiplatform
  environment for analyzing dynamical systems,'' {\em User’s Guide, Dynamical
  Systems Laboratory, CWI, Amsterdam, Netherlands, Available by anonymous ftp
  from ftp.cwi.nl/pub/CONTENT}, 1995-1997.

\bibitem{dhooge2003}
A.~Dhooge, W.~Govaerts, and Y.~A. Kuznetsov, ``Matcont: a matlab package for
  numerical bifurcation analysis of odes,'' {\em ACM Transactions on
  Mathematical Software (TOMS)}, vol.~29, no.~2, pp.~141--164, 2003.

\bibitem{dankowicz2011}
H.~Dankowicz and F.~Schilder, ``An extended continuation problem for
  bifurcation analysis in the presence of constraints,'' {\em Journal of
  Computational and Nonlinear Dynamics}, vol.~6, no.~3, 2011.

\bibitem{krylov1943}
N.~M. Krylov and N.~N. Bogoliubov, {\em Introduction to non-linear mechanics}.
\newblock No.~11, Princeton University Press, 1943.

\bibitem{urabe1965}
M.~Urabe, ``Galerkin's procedure for nonlinear periodic systems,'' {\em Archive
  for Rational Mechanics and Analysis}, vol.~20, no.~2, pp.~120--152, 1965.

\bibitem{karkar2014}
S.~Karkar, B.~Cochelin, and C.~Vergez, ``A comparative study of the harmonic
  balance method and the orthogonal collocation method on stiff nonlinear
  systems,'' {\em Journal of Sound and Vibration}, vol.~333, no.~12,
  pp.~2554--2567, 2014.

\bibitem{lau1981}
S.~Lau and Y.~Cheung, ``Amplitude incremental variational principle for
  nonlinear vibration of elastic systems,'' {\em Journal of Applied Mechanics},
  vol.~48, no.~4, pp.~959--964, 1981.

\bibitem{cheung1982}
Y.~Cheung and S.~Lau, ``Incremental time—space finite strip method for
  non-linear structural vibrations,'' {\em Earthquake Engineering \& Structural
  Dynamics}, vol.~10, no.~2, pp.~239--253, 1982.

\bibitem{sze2005}
K.~Sze, S.~Chen, and J.~Huang, ``The incremental harmonic balance method for
  nonlinear vibration of axially moving beams,'' {\em Journal of Sound and
  Vibration}, vol.~281, no.~3, pp.~611--626, 2005.

\bibitem{cameron1989}
T.~Cameron and J.~Griffin, ``An alternating frequency/time domain method for
  calculating the steady-state response of nonlinear dynamic systems,'' {\em
  Journal of Applied Mechanics}, vol.~56, no.~1, pp.~149--154, 1989.

\bibitem{cardona1998}
A.~Cardona, A.~Lerusse, and M.~G{\'e}radin, ``Fast fourier nonlinear vibration
  analysis,'' {\em Computational Mechanics}, vol.~22, no.~2, pp.~128--142,
  1998.

\bibitem{sinou2007}
J.-J. Sinou and A.~Lees, ``A non-linear study of a cracked rotor,'' {\em
  European Journal of Mechanics-A/Solids}, vol.~26, no.~1, pp.~152--170, 2007.

\bibitem{jaumouille2010}
V.~Jaumouill{\'e}, J.-J. Sinou, and B.~Petitjean, ``An adaptive harmonic
  balance method for predicting the nonlinear dynamic responses of mechanical
  systems - {Application} to bolted structures,'' {\em Journal of Sound and
  Vibration}, vol.~329, no.~19, pp.~4048--4067, 2010.

\bibitem{guillen1999}
J.~Guillen and C.~Pierre, ``An efficient, hybrid, frequency-time domain method
  for the dynamics of large-scale dry-friction damped structural systems.,'' in
  {\em IUTAM Symposium on Unilateral Multibody Contacts}, pp.~169--178,
  Springer, 1999.

\bibitem{poudou2003}
O.~Poudou and C.~Pierre, ``Hybrid frequency-time domain methods for the
  analysis of complex structural systems with dry friction damping,'' in {\em
  Collection of Technical Papers—AIAA/ASME/ASCE/AHS/ASC Structures,
  Structural Dynamics and Materials Conference}, pp.~111--124, 2003.

\bibitem{vongroll2001}
G.~von Groll and D.~J. Ewins, ``The harmonic balance method with arc-length
  continuation in rotor/stator contact problems,'' {\em Journal of Sound and
  Vibration}, vol.~241, no.~2, pp.~223--233, 2001.

\bibitem{arquier2007}
R.~Arquier, {\em Une m{\'e}thode de calcul des modes de vibrations
  non-lin{\'e}aires de structures}.
\newblock PhD thesis, Universit{\'e} de la m{\'e}diterran{\'e}e (Aix-Marseille
  II), Marseille, France, 2007.

\bibitem{cochelin2009}
B.~Cochelin and C.~Vergez, ``A high order purely frequency-based harmonic
  balance formulation for continuation of periodic solutions,'' {\em Journal of
  sound and vibration}, vol.~324, no.~1, pp.~243--262, 2009.

\bibitem{grolet2012}
A.~Grolet and F.~Thouverez, ``On a new harmonic selection technique for
  harmonic balance method,'' {\em Mechanical Systems and Signal Processing},
  vol.~30, pp.~43--60, 2012.

\bibitem{kundert1986}
K.~S. Kundert and A.~Sangiovanni-Vincentelli, ``Simulation of nonlinear
  circuits in the frequency domain,'' {\em IEEE Transactions on Computer-Aided
  Design of Integrated Circuits and Systems}, vol.~5, no.~4, pp.~521--535,
  1986.

\bibitem{genesio1992}
R.~Genesio and A.~Tesi, ``Harmonic balance methods for the analysis of chaotic
  dynamics in nonlinear systems,'' {\em Automatica}, vol.~28, no.~3,
  pp.~531--548, 1992.

\bibitem{stanton2012}
S.~C. Stanton, B.~A. Owens, and B.~P. Mann, ``Harmonic balance analysis of the
  bistable piezoelectric inertial generator,'' {\em Journal of Sound and
  Vibration}, vol.~331, no.~15, pp.~3617--3627, 2012.

\bibitem{fang2012}
C.-C. Fang, ``Critical conditions for a class of switched linear systems based
  on harmonic balance: applications to dc-dc converters,'' {\em Nonlinear
  Dynamics}, vol.~70, no.~3, pp.~1767--1789, 2012.

\bibitem{shen1959}
S.~F. Shen, ``An approximate analysis of nonlinear flutter problems,'' {\em
  Journal of the Aerospace Sciences}, vol.~26, no.~1, pp.~25--32, 1959.

\bibitem{liu2005}
L.~Liu and E.~H. Dowell, ``Harmonic balance approach for an airfoil with a
  freeplay control surface,'' {\em AIAA journal}, vol.~43, no.~4, pp.~802--815,
  2005.

\bibitem{lee2005}
B.~Lee, L.~Liu, and K.~Chung, ``Airfoil motion in subsonic flow with strong
  cubic nonlinear restoring forces,'' {\em Journal of Sound and Vibration},
  vol.~281, no.~3, pp.~699--717, 2005.

\bibitem{dimitriadis2008}
G.~Dimitriadis, ``Continuation of higher-order harmonic balance solutions for
  nonlinear aeroelastic systems,'' {\em Journal of aircraft}, vol.~45, no.~2,
  pp.~523--537, 2008.

\bibitem{hall2002}
K.~C. Hall, J.~P. Thomas, and W.~S. Clark, ``Computation of unsteady nonlinear
  flows in cascades using a harmonic balance technique,'' {\em AIAA journal},
  vol.~40, no.~5, pp.~879--886, 2002.

\bibitem{gopinath2007}
A.~Gopinath, E.~Van Der~Weide, J.~J. Alonso, A.~Jameson, K.~Ekici, and K.~C.
  Hall, ``Three-dimensional unsteady multi-stage turbomachinery simulations
  using the harmonic balance technique,'' in {\em 45th AIAA Aerospace Sciences
  Meeting and Exhibit}, no.~892, 2007.

\bibitem{clark2015}
S.~T. Clark, F.~M. Besem, R.~E. Kielb, and J.~P. Thomas, ``Developing a
  reduced-order model of nonsynchronous vibration in turbomachinery using
  proper-orthogonal decomposition methods,'' {\em Journal of Engineering for
  Gas Turbines and Power}, vol.~137, no.~5, p.~052501, 2015.

\bibitem{ekici2008}
K.~Ekici, K.~C. Hall, and E.~H. Dowell, ``Computationally fast harmonic balance
  methods for unsteady aerodynamic predictions of helicopter rotors,'' {\em
  Journal of Computational Physics}, vol.~227, no.~12, pp.~6206--6225, 2008.

\bibitem{cardona1994}
A.~Cardona, T.~Coune, A.~Lerusse, and M.~Geradin, ``A multiharmonic method for
  non-linear vibration analysis,'' {\em International Journal for Numerical
  Methods in Engineering}, vol.~37, no.~9, pp.~1593--1608, 1994.

\bibitem{petrov2003}
E.~Petrov and D.~Ewins, ``Analytical formulation of friction interface elements
  for analysis of nonlinear multi-harmonic vibrations of bladed disks,'' {\em
  Journal of turbomachinery}, vol.~125, no.~2, pp.~364--371, 2003.

\bibitem{suss2015}
D.~S{\"u}{\ss} and K.~Willner, ``Investigation of a jointed friction oscillator
  using the multiharmonic balance method,'' {\em Mechanical Systems and Signal
  Processing}, vol.~52, pp.~73--87, 2015.

\bibitem{woo2000}
K.-C. Woo, A.~A. Rodger, R.~D. Neilson, and M.~Wiercigroch, ``Application of
  the harmonic balance method to ground moling machines operating in periodic
  regimes,'' {\em Chaos, Solitons \& Fractals}, vol.~11, no.~15,
  pp.~2515--2525, 2000.

\bibitem{peter2014}
S.~Peter, P.~Reuss, and L.~Gaul, ``Identification of sub-and higher harmonic
  vibrations in vibro-impact systems,'' in {\em Nonlinear Dynamics, Volume 2},
  pp.~131--140, Springer, 2014.

\bibitem{lewandowski1997}
R.~Lewandowski, ``Computational formulation for periodic vibration of
  geometrically nonlinear structures—part 2: numerical strategy and
  examples,'' {\em International journal of solids and structures}, vol.~34,
  no.~15, pp.~1949--1964, 1997.

\bibitem{ribeiro1999}
P.~Ribeiro and M.~Petyt, ``Geometrical non-linear, steady state, forced,
  periodic vibration of plates, part i: model and convergence studies,'' {\em
  Journal of Sound and Vibration}, vol.~226, no.~5, pp.~955--983, 1999.

\bibitem{barillon2013}
F.~Barillon, J.-J. Sinou, J.-M. Duffal, and L.~J{\'e}z{\'e}quel, ``Non--linear
  dynamics of a whole vehicle finite element model using a harmonic balance
  method,'' {\em International Journal of Vehicle Design}, vol.~63, no.~4,
  pp.~387--403, 2013.

\bibitem{senechal2014}
A.~S{\'e}n{\'e}chal, B.~Petitjean, and L.~Zoghai, ``Development of a numerical
  tool for industrial structures with local nonlinearities,'' in {\em
  Proceedings of the 26th International Conference on Noise and Vibration
  engineering (ISMA2014)}, (Leuven, Belgium), 2014.

\bibitem{claeys2014}
M.~Claeys, J.-J. Sinou, J.-P. Lambelin, and B.~Alcoverro, ``Multi-harmonic
  measurements and numerical simulations of nonlinear vibrations of a beam with
  non-ideal boundary conditions,'' {\em Communications in Nonlinear Science and
  Numerical Simulation}, 2014.

\bibitem{claeys2014b}
M.~Claeys, J.-J. Sinou, J.-P. Lambelin, and R.~Todeschini, ``Experimental and
  numerical study of the nonlinear vibrations of an assembly with friction
  joints,'' in {\em Proceedings of the 26th International Conference on Noise
  and Vibration engineering (ISMA2014)}, (Leuven, Belgium), 2014.

\bibitem{kerschen2009}
G.~Kerschen, M.~Peeters, J.-C. Golinval, and A.~F. Vakakis, ``Nonlinear normal
  modes, part i: A useful framework for the structural dynamicist,'' {\em
  Mechanical Systems and Signal Processing}, vol.~23, no.~1, pp.~170--194,
  2009.

\bibitem{krack2013b}
M.~Krack, L.~Panning-von Scheidt, and J.~Wallaschek, ``A high-order harmonic
  balance method for systems with distinct states,'' {\em Journal of Sound and
  Vibration}, vol.~332, no.~21, pp.~5476--5488, 2013.

\bibitem{detroux2014}
T.~Detroux, L.~Renson, and G.~Kerschen, ``The harmonic balance method for
  advanced analysis and design of nonlinear mechanical systems,'' in {\em
  Proceedings of the 32th International Modal Analysis Conference (IMAC)},
  (Orlando, FL, USA), 2014.

\bibitem{kuether2014}
R.~J. Kuether, M.~R. Brake, and M.~S. Allen, ``Evaluating convergence of
  reduced order models using nonlinear normal modes,'' in {\em Model Validation
  and Uncertainty Quantification, Volume 3}, pp.~287--300, Springer, 2014.

\bibitem{laxalde2009}
D.~Laxalde and F.~Thouverez, ``Complex non-linear modal analysis for mechanical
  systems: Application to turbomachinery bladings with friction interfaces,''
  {\em Journal of sound and vibration}, vol.~322, no.~4, pp.~1009--1025, 2009.

\bibitem{krack2013}
M.~Krack, L.~Panning-von Scheidt, J.~Wallaschek, C.~Siewert, and A.~Hartung,
  ``Reduced order modeling based on complex nonlinear modal analysis and its
  application to bladed disks with shroud contact,'' {\em Journal of
  Engineering for Gas Turbines and Power}, vol.~135, no.~10, p.~102502, 2013.

\bibitem{schilder2005}
F.~Schilder, H.~M. Osinga, and W.~Vogt, ``Continuation of quasi-periodic
  invariant tori,'' {\em SIAM Journal on Applied Dynamical Systems}, vol.~4,
  no.~3, pp.~459--488, 2005.

\bibitem{peletan2014}
L.~Peletan, S.~Baguet, M.~Torkhani, and G.~Jacquet-Richardet, ``Quasi-periodic
  harmonic balance method for rubbing self-induced vibrations in rotor--stator
  dynamics,'' {\em Nonlinear Dynamics}, vol.~78, no.~4, pp.~2501--2515, 2014.

\bibitem{guskov2012}
M.~Guskov and F.~Thouverez, ``Harmonic balance-based approach for
  quasi-periodic motions and stability analysis,'' {\em Journal of Vibration
  and Acoustics}, vol.~134, no.~3, p.~031003, 2012.

\bibitem{liao2014}
H.~Liao, ``Global resonance optimization analysis of nonlinear mechanical
  systems: Application to the uncertainty quantification problems in rotor
  dynamics,'' {\em Communications in Nonlinear Science and Numerical
  Simulation}, vol.~19, no.~9, pp.~3323--3345, 2014.

\bibitem{grolet2013}
A.~Grolet and F.~Thouverez, ``Vibration of mechanical systems with geometric
  nonlinearities: Solving harmonic balance equations with groebner basis and
  continuations methods,'' in {\em Proceedings of the Colloquium Calcul des
  Structures et Modélisation CSMA}, (Giens, France), 2013.

\bibitem{hill1886}
G.~W. Hill, ``On the part of the motion of the lunar perigee which is a
  function of the mean motions of the sun and moon,'' {\em Acta mathematica},
  vol.~8, no.~1, pp.~1--36, 1886.

\bibitem{lanza2007}
V.~Lanza, M.~Bonnin, and M.~Gilli, ``On the application of the describing
  function technique to the bifurcation analysis of nonlinear systems,'' {\em
  Circuits and Systems II: Express Briefs, IEEE Transactions on}, vol.~54,
  no.~4, pp.~343--347, 2007.

\bibitem{bonani1999}
F.~Bonani and M.~Gilli, ``Analysis of stability and bifurcations of limit
  cycles in chua's circuit through the harmonic-balance approach,'' {\em
  Circuits and Systems I: Fundamental Theory and Applications, IEEE
  Transactions on}, vol.~46, no.~8, pp.~881--890, 1999.

\bibitem{traversa2008}
F.~Traversa, F.~Bonani, and S.~D. Guerrieri, ``A frequency-domain approach to
  the analysis of stability and bifurcations in nonlinear systems described by
  differential-algebraic equations,'' {\em International Journal of Circuit
  Theory and Applications}, vol.~36, no.~4, pp.~421--439, 2008.

\bibitem{lazarus2010}
A.~Lazarus and O.~Thomas, ``A harmonic-based method for computing the stability
  of periodic solutions of dynamical systems,'' {\em Comptes Rendus
  M{\'e}canique}, vol.~338, no.~9, pp.~510--517, 2010.

\bibitem{traversa2012}
F.~L. Traversa and F.~Bonani, ``Improved harmonic balance implementation of
  floquet analysis for nonlinear circuit simulation,'' {\em AEU-International
  Journal of Electronics and Communications}, vol.~66, no.~5, pp.~357--363,
  2012.

\bibitem{piccardi1994}
C.~Piccardi, ``Bifurcations of limit cycles in periodically forced nonlinear
  systems: The harmonic balance approach,'' {\em Circuits and Systems I:
  Fundamental Theory and Applications, IEEE Transactions on}, vol.~41, no.~4,
  pp.~315--320, 1994.

\bibitem{piccardi1996}
C.~Piccardi, ``Harmonic balance analysis of codimension-2 bifurcations in
  periodic systems,'' {\em IEEE Transactions on Circuits and Systems -- Part 1:
  Fundamental Theory and Applications}, vol.~43, pp.~1015--1018, 1996.

\bibitem{powell1970}
M.~J. Powell, ``A hybrid method for nonlinear equations,'' {\em Numerical
  methods for nonlinear algebraic equations}, vol.~7, pp.~87--114, 1970.

\bibitem{nacivet2003}
S.~Nacivet, C.~Pierre, F.~Thouverez, and L.~Jezequel, ``A dynamic lagrangian
  frequency--time method for the vibration of dry-friction-damped systems,''
  {\em Journal of Sound and Vibration}, vol.~265, no.~1, pp.~201--219, 2003.

\bibitem{Hwang1991}
J.~L. Hwang and T.~N. Shiau, ``An application of the generalized polynomial
  expansion method to nonlinear rotor bearing systems,'' {\em Journal of
  Vibration and Acoustics}, vol.~113, no.~3, pp.~299--308, 1991.

\bibitem{Xie1996}
G.~Xie and J.~Y. Lou, ``Alternating frequency/coefficient (afc) technique in
  the trigonometric collocation method,'' {\em International journal of
  non-linear mechanics}, vol.~31, no.~4, pp.~531--545, 1996.

\bibitem{narayanan1998}
S.~Narayanan and P.~Sekar, ``A frequency domain based numeric--analytical
  method for non-linear dynamical systems,'' {\em Journal of sound and
  vibration}, vol.~211, no.~3, pp.~409--424, 1998.

\bibitem{duan2005}
C.~Duan and R.~Singh, ``Super-harmonics in a torsional system with dry friction
  path subject to harmonic excitation under a mean torque,'' {\em Journal of
  Sound and Vibration}, vol.~285, no.~4, pp.~803--834, 2005.

\bibitem{kim2005}
T.~Kim, T.~Rook, and R.~Singh, ``Super-and sub-harmonic response calculations
  for a torsional system with clearance nonlinearity using the harmonic balance
  method,'' {\em Journal of Sound and Vibration}, vol.~281, no.~3,
  pp.~965--993, 2005.

\bibitem{peletan2013}
L.~Peletan, S.~Baguet, M.~Torkhani, and G.~Jacquet-Richardet, ``A comparison of
  stability computational methods for periodic solution of nonlinear problems
  with application to rotordynamics,'' {\em Nonlinear Dynamics}, vol.~72,
  no.~3, pp.~671--682, 2013.

\bibitem{moore2005}
G.~Moore, ``Floquet theory as a computational tool,'' {\em SIAM journal on
  numerical analysis}, vol.~42, no.~6, pp.~2522--2568, 2005.

\bibitem{govaerts2011}
W.~Govaerts, Y.~A. Kuznetsov, V.~De~Witte, A.~Dhooge, H.~Meijer, W.~Mestrom,
  A.~Riet, and B.~Sautois, ``Matcont and cl matcont: Continuation toolboxes in
  matlab,'' {\em Gent University and Utrech University, Tech. Rep}, 2011.

\bibitem{seydel2010}
R.~Seydel, {\em Practical bifurcation and stability analysis}.
\newblock Springer, 2010.

\bibitem{guckenheimer1997}
J.~Guckenheimer, M.~Myers, and B.~Sturmfels, ``Computing hopf bifurcations i,''
  {\em SIAM Journal on Numerical Analysis}, vol.~34, no.~1, pp.~1--21, 1997.

\bibitem{beyn2002}
W.-J. Beyn, A.~Champneys, E.~Doedel, W.~Govaerts, Y.~A. Kuznetsov, and
  B.~Sandstede, ``Numerical continuation, and computation of normal forms,''
  {\em Handbook of dynamical systems}, vol.~2, pp.~149--219, 2002.

\bibitem{vanderaa2007}
N.~Van Der~Aa, H.~Ter~Morsche, and R.~Mattheij, ``Computation of eigenvalue and
  eigenvector derivatives for a general complex-valued eigensystem,'' {\em
  Electronic Journal of Linear Algebra}, vol.~16, no.~1, pp.~300--314, 2007.

\bibitem{noel2014}
J.~No{\"e}l, L.~Renson, and G.~Kerschen, ``Complex dynamics of a nonlinear
  aerospace structure: Experimental identification and modal interactions,''
  {\em Journal of Sound and Vibration}, vol.~333, no.~12, pp.~2588--2607, 2014.

\bibitem{renson2014}
L.~Renson, J.~No{\"e}l, and G.~Kerschen, ``Complex dynamics of a nonlinear
  aerospace structure: numerical continuation and normal modes,'' {\em
  Nonlinear Dynamics}, vol.~79, no.~2, pp.~1293--1309, 2015.

\end{thebibliography}
\bibliographystyle{ieeetr}

\end{document}